\documentclass[a4paper,twoside]{article}
\usepackage{mathrsfs}
\usepackage{amsfonts}
\usepackage{amsfonts,amsmath,amssymb,amsthm}
\usepackage{latexsym,amscd,graphics,color}
\usepackage{amsbsy}

\def\Rep{{\mbox{\bf Rep}}}
\def\Mod{\mbox{\bf Mod}}
\def\Hom{\mbox{Hom}}
\def\Tr{\mbox{Tr}}

\def\rad{\mbox{rad}}
\def\End{\mbox{End}}
\def\dim{\mbox{dim}}
\def\Aut{\mbox{Aut}}
\def\id{\mbox{id}}
\def\Im{\mbox{Im}}
\def\Ker{\mbox{Ker}}
\def\irr{\mbox{irr}}

\begin{document}
\baselineskip=16pt
\title{Skew group algebras of deformed preprojective algebras
\footnotetext {$\ast$ E-mail: houbo@emails.bjut.edu.cn; ~~
$\dag$ E-mail: slyang@bjut.edu.cn }
\footnotetext { $\dag$ Supported by National Natural Science
Foundation of China (Grants. 10771014) } }

\author{\large Bo Hou$^{\ast}$, ~Shilin Yang$^{\dag}$  }

\date{\footnotesize College of Applied Science, Beijing University of Technology,
 Beijing 100124, China }

\maketitle

\begin{center}

\begin{minipage}{14cm}

{\bf Abstract}: Suppose that $Q$ is a finite quiver and $G\subseteq \Aut(Q)$
is a finite group, $k$ is an algebraic closed field whose
characteristic does not divide the order of $G$. For any algebra
$\Lambda=kQ/{\mathcal {I}}$, $\mathcal {I}$ is an arbitrary ideal of
path algebra $kQ$, we give all the indecomposable $\Lambda
G$-modules from indecomposable $\Lambda$-modules when $G$ is
abelian. In particular, we apply this result to the deformed
preprojective algebra $\Pi_{Q}^{\lambda}$, and get a reflection
functor for the module category of $\Pi_{Q}^{\lambda}G$.
Furthermore, we construct a new quiver $Q_{G}$ and prove that
$\Pi_{Q}^{\lambda}G$ is Morita equivalent to $\Pi_{Q_{G}}^{\eta}$
for some $\eta$.



\end{minipage}

\end{center}

\section{Introduction}
Let $k$ be a field and $Q=(I, A)$ a finite quiver with a finite
set $I$ of vertices  and a finite set $A$ of arrows. For each arrow
$a\in A$, its head vertex and tail vertex are denoted by $h(a)$ and
$t(a)$ respectively. For each $i\in I$ we denote by $e_{i}$ the primitive
idempotent of path algebra $kQ$ corresponding to vertex $i$.
 The double quiver of $Q$ is denoted by $\overline{Q}=(I, \overline{A})$,
 which
is obtained from $Q$ by adding an arrow $\alpha': j\rightarrow i$
whenever there is an arrow $\alpha : i\rightarrow j$ in $Q$.
Given  a
$\lambda=(\lambda_{i})_{i\in I}\in k^{I}$, let
$\mathcal {I}$ be the ideal of path algebra $k\overline{Q}$
generated by the deformed preprojective relations
$$\sum_{a\in A, ~h(a)=i}aa'-\sum_{a\in A, ~t(a)
=i}a'a-\lambda_{i}e_{i}, \qquad \forall ~i\in I,$$ The quotient algebra
$\Pi_{Q}^{\lambda}:=k\overline{Q}/{\mathcal {I}}$ is called a
deformed preprojective algebra of $Q$. In particular,  if
$\lambda=0$,
$\Pi_{Q}^{\lambda}$ is called a preprojective algebra of $Q$,
which is denoted by $\Pi_{Q}$.

Preprojective algebras were introduced by Gelfand and Ponomarev
\cite{GP} to study  preprojective representations of a finite quiver
without oriented cycle. Subsequently it was discovered that
preprojective algebras occur naturally in very diverse contexts. For
example, they have been used by Kronheimer \cite{K} to deal with the
problems in differential geometry, also by Lusztig \cite{L1,L2} in
his perverse sheaf approach to quantum groups. In 1998, deformed
preprojective algebras were introduced by Crawley-Boevey and Holland
\cite{CH} to study noncommutative deformations of Kleinian
singularities. They are generalized preprojective algebras. In the
papers \cite{Y}, \cite{HY}, the quotient algebras of deformed
preprojective algebras were used to realize the restricted quantum
group $U_{q}({\mathfrak{sl}}_{2})$. Recently, Demonet proved that
skew  group algebras of preprojective algebras are Morita equivalent
to preprojective algebras \cite{De}.

\medskip
Skew group algebras appear in connection with the study of
singularities \cite{A, AR}. They are studied by many authors (see
for example \cite{FR}, \cite{RR}, \cite{MR},\cite{ALR}).
For the study of the
representation theory of skew group algebras, the
readers refer to \cite{RR} (see  also \cite{ARS}). Let
$\Lambda$ be a $k$-algebra with a group  $G$ acting on $\Lambda$, then the
skew group algebra $\Lambda G$ is an associative $k$-algebra whose
underlying $k$-vector space is spanned by the elements $ga$ for
$g\in G$ and $a\in \Lambda$, and the multiplication is defined by
$$ga\cdot hb=ghh^{-1}(a)b,$$ for all $g, h \in G$ and $a, b \in
\Lambda$ (see \cite{RR}).

 It is well-known that skew group algebras
$\Lambda G$ retain many features from $\Lambda$, such as
representation-finite,  an Auslander algebra,  a piecewise
hereditary algebra, a Nakayama algebra, or  a Koszul
algebra, and so on (see \cite{ARS, FR, MR, ALR, DLS}). But there are many
problems unsolved  about the relationships between $\Lambda$ and $\Lambda G$.
The aim of this paper is to consider
what kind of $\Lambda$-modules can induce $\Lambda G$-modules,
and if a $\Lambda$-module induce a $\Lambda G$-module, how many
isomorphism classes of such induced $\Lambda G$-module?

As far as we know, in \cite{ZL}  the authors studied the relationships between
$\Lambda G$ and $\Lambda$
under the assumption that $\Lambda$ is a path algebra and $G$ is
a cyclic group.
In Section 2, we will
discuss the same problems  for an arbitrary algebra
$\Lambda=kQ/{\mathcal {I}}$ and an arbitrary finite group
$G\subseteq \Aut(Q,{\mathcal {I}})$, where $\Aut(Q,{\mathcal
{I}})=\{\alpha\in \Aut(Q)\mid \alpha({\mathcal {I}})={\mathcal
{I}}\}$.  $\Aut(Q,{\mathcal {I}})$  is called the automorphism
group of bound quiver
$(Q,{\mathcal {I}})$. If $k$ is an algebraic closed field whose
characteristic does not divide the order of $G$, it is shown that
$\Lambda$-module $Y$ is a $\Lambda G$-module if and only if $Y$ is a
$G$-invariant $\Lambda$-module. Sequentially, for any indecomposable
$\Lambda$-module $X$, if $H_{X}=\{g\in G \mid {^{g}X}\cong X\}$ is
an abelian group,  we describe all the $\Lambda$-module structures
on the $G$-invariant $\Lambda$-module $\bigoplus_{g\in
G_{X}}{^{g}X}$ by the irreducible decomposition
$kH_{X}=\bigoplus_{i=1}^{r}\rho_{i}$ as $H_{X}$-representations.

\medskip
{\bf Theorem 1.1.} Suppose $X$ is a finite dimensional
indecomposable $\Lambda$-module, $H_{X}$ is abelian. Then for any
$\Lambda G$-module $Y$, if $Y\cong\bigoplus_{g\in G_{X}}{^{g}X}$ as
$\Lambda$-modules, there exists a unique $i\in \{1,2,\cdots,r\}$
such that $Y\cong \Lambda G\otimes_{\Lambda
H_{X}}(L_{i}\otimes_{k}X)$. That is, there are $r$ non-isomorphic
$\Lambda G$-modules induced from the indecomposable $G$-invariant
$\Lambda$-module $\bigoplus_{g\in G_{X}}{^{g}X}$.

\medskip
In fact, if $G\subseteq \Aut(Q,{\mathcal {I}})$ is a finite abelian
group, we can get all the indecomposable $\Lambda G$-modules from
indecomposable $\Lambda$-modules.

\medskip
{\bf Theorem 1.2.} If a finite group $G\subseteq \Aut(Q, {\mathcal
{I}})$ is abelian, then a finite dimensional $\Lambda$-module $Y$ is
an indecomposable $\Lambda G$-module if and only if $Y$ is an
indecomposable $G$-invariant $\Lambda$-module, that is,
$Y\cong\bigoplus_{g\in G_{X}}{^{g}X}$, for some indecomposable
$\Lambda$-module $X$.

\medskip
According to Theorem 1.1 and Theorem 1.2, for a skew group algebra
$\Lambda G$ with $G$ is abelian, all finite dimensional $\Lambda
G$-modules can be obtained from indecomposable $\Lambda$-modules.
The number of non-isomorphic indecomposable $\Lambda G$-modules
induced from the same indecomposable $G$-invariant $\Lambda$-module
can be given. In particular, all finite dimensional
$\Pi_{Q}^{\lambda}G$-modules from indecomposable
$\Pi_{Q}^{\lambda}$-modules can be constructed.

\medskip
In \cite{CH}, Crawley-Boevey and Holland introduced an interesting
reflection functor for deformed preprojective algebras. In Section 3,
we will define a
reflection $r_{{\mathcal {O}}_{i}}$ on dimension vectors and a
reflection $d_{{\mathcal {O}}_{i}}$ on $k^{I}$. Then we use
Crawley-Boevey and Holland$^\prime$s reflection functor to construct a
reflection functor for skew group algebras of deformed preprojective
algebras.

\medskip
{\bf Theorem 1.3.} Suppose $\lambda=(\lambda_{i})_{i\in I}$
satisfies $\lambda_{i}=\lambda_{j}$ for any $i, j\in I$ in the same
$G$-orbit, $i$ is a loop-free vertex and $\lambda_{i}\neq 0$. If $G$
is abelian and the action of $G$ is admissible, then there is an
equivalence
$$F_{i}:\Pi_{Q}^{\lambda}G\mbox{-{\bf mod}}\longrightarrow
\Pi_{Q}^{d_{{\mathcal {O}}_{i}}\lambda}G\mbox{-{\bf mod}},$$ that
acts as the reflection $r_{{\mathcal {O}}_{i}}$ on dimension
vectors.

\medskip
Recently, Demonet has introduced the definition of group species and
their representations, and proved that skew group algebras of
preprojective algebras are Morita equivalent to preprojective
algebras \cite{De}. In fact, it also holds for deformed preprojective
algebras. In Section 4, we will give the group species corresponding
to quiver $Q=(I, A)$ with $G\subseteq \Aut(Q)$ and construct a new
quiver $Q_{G}=\big(I_{G}, ~(A'_{(i, \rho), (j, \sigma)})_{((i,
\rho), (j, \sigma))\in I^{2}_{G}}\big)$. We get the following
Morita equivalence.

\medskip
{\bf Theorem 1.4.} If $\lambda=(\lambda_{i})\in k^{I}$ satisfies
$\lambda_{i}=\lambda_{j}$ for any $i, j\in I$ in the same $G$-orbit,
then there is an equivalence of categories
$$\Pi_{Q}^{\lambda}G\mbox{-{\bf Mod}}\simeq
\Pi_{Q_{G}}^{\eta}\mbox{-{\bf Mod}},$$ where $\eta=(\eta_{(i,
\rho)})\in k^{I_{G}}$, and $\eta_{(i, \rho)}=\frac{\dim_{k}
\rho}{|G_{i}|} \lambda_{i}$.
In particular, $\Pi_{Q}G\mbox{-{\bf Mod}}\simeq
\Pi_{Q_{G}}\mbox{-{\bf Mod}}$.

\medskip
Consequently, the study of the skew group algebras of deformed
preprojective algebras can be reduced to the study of deformed
preprojective algebras. Finally, as an example, the quiver
$Q_{G}$ is given for any Dynkin quiver $Q$ with $G\subseteq \Aut(Q)$.

\medskip
In this paper, we suppose that  $Q=(I, A)$ is a finite quiver and
$G\subseteq \Aut(Q)$ is a finite group with unit $e$. Let  $k$ be an
algebraic closed field whose characteristic does not divide the
order of $G$, i.e., char$k\nmid |G|$. For any algebra $\Lambda$, we
denote by $\Lambda$-${\bf Mod}$ the category of $\Lambda$-modules,
by $\Lambda$-${\bf mod}$ the full subcategory of $\Lambda$-${\bf
Mod}$ consisting of finite dimensional $\Lambda$-modules. The readers refer to
\cite{ARS} for the knowledge of representation of quiver and bound quiver.

\section{Module category of skew group algebras}
In this section, we suppose that $\mathcal {I}$ is an arbitrary
ideal of path algebra $kQ$ ( $\mathcal {I}$ is not necessarily
admissible), $\Lambda=kQ/{\mathcal {I}}$, and $G\subseteq \Aut(Q,
{\mathcal {I}})$ is a finite group. All modules (or representations)
are finite dimensional. ${\bf Rep}(Q, {\mathcal {I}})$ is the
category of finite dimensional representations of bound quiver $(Q,
{\mathcal {I}})$.

\medskip
Let $X$ be a $\Lambda$-module, $g\in G$. We define a twisted
$\Lambda$-module $^{g}X$ on $X$ as follows: as a $k$-vector space
$^gX=X$,  the action on $^gX$ is given by $a\cdot x=g^{-1}(a)x$ for
all $a\in \Lambda$. Let $\varphi:X\rightarrow Y$ be a
$\Lambda$-module homomorphism, then the map $^{g}\varphi=\varphi$
can be viewed as the twisted $\Lambda$-module homomorphism:
${^{g}X}\to {^{g}Y}$. Indeed, for $x\in X$ and $a\in \Lambda$, we
have $\varphi(a\cdot
x)=\varphi(g^{-1}(a)x)=g^{-1}(a)\varphi(x)=a\cdot\varphi(x)$.

It is well-known that categories ${\bf Rep}(Q, {\mathcal {I}})$ and
$\Lambda$-{\bf mod} are equivalent. Therefore, for any $X\in
\Lambda$-{\bf mod}, there is one and only one representation
$(X_{i}, ~X_{a})_{i\in I, a\in A}$ of bound quiver $(Q, {\mathcal
{I}})$ corresponding to $X$. It is easy to see that
$(^{g}X_{i}=X_{g^{-1}(i)},~~ {^{g}X}_{a}=X_{g^{-1}(a)})_{i\in I,
a\in A}$ is the representation of bound quiver $(Q, {\mathcal {I}})$
corresponding to $^{g}X$. For  $X, Y\in \Lambda$-{\bf mod} and
homomorphism $\varphi:X\rightarrow Y$, we define a functor $F_g$ by
$F_{g}(X)={^{g}X}$ and $F_{g}(\varphi)={^{g}\varphi}$.  One can
check that
$$F_{g} : \Lambda\mbox{-{\bf mod}}\longrightarrow
\Lambda\mbox{-{\bf mod}}$$ is an equivalence functor. The inverse is
$F_{g^{-1}}$. It follows that for each $g\in G$, $X\in \Lambda$-{\bf mod} is
indecomposable (projective, injective, simple) if and only if
so is $^{g}X$.

\medskip
{\bf Definition 2.1.} (1) Given any $g\in G$, a representation $X$
of bound quiver $(Q, {\mathcal {I}})$ (or $X\in \Lambda$-{\bf mod})
 is said to be {\bf $g$-invariant} if $^{g}X\cong X$ as
representations of bound quiver $(Q, {\mathcal {I}})$.

(2) A representation $X$ of bound quiver $(Q, {\mathcal {I}})$ is
said to be {\bf $G$-invariant} if $X$ is $g$-invariant for any $g\in G$.

(3) A $G$-invariant representation $X$ of bound quiver $(Q,
{\mathcal {I}})$ is called {\bf indecomposable $G$-invariant} if $X$
is nonzero and $X$ can not be decomposed as $X\cong X'\oplus X''$,
where $X'$ and $X''$ are nonzero $G$-invariant representations.

\medskip
For each $X\in \Lambda$-{\bf mod}, let $H_{X}=\{g\in G \mid
{^{g}X}\cong X\}$. Clearly, $H_{X}$ is a subgroup of $G$. If we
denote by $G_{X}$ a complete set of left coset representatives of
$H_{X}$ in $G$, then we have

\medskip
{\bf Lemma 2.1.} Any indecomposable $G$-invariant representation of
bound quiver $(Q, {\mathcal {I}})$ is
of the form
$\bigoplus_{g\in G_{X}}{^{g}X}$ where $X$ is an indecomposable
$(Q, {\mathcal {I}})$-representation.
Moreover, the Krull-Schmidt theorem holds for $G$-invariant
representations.

\medskip
{\bf Proof.} Let $X$ be an indecomposable $(Q, {\mathcal
{I}})$-representation. Since $\bigoplus_{g\in E}{^{g}X}$ is not
$G$-invariant for any proper subset $E$ of $G_{X}$ and
$\bigoplus_{g\in G_{X}}{^{g}X}\cong \bigoplus_{g\in G_{X}}{^{hg}X}$
for any $h\in G$,  $\bigoplus_{g\in G_{X}}{^{g}X}$ is an
indecomposable $G$-invariant $(Q, {\mathcal {I}})$-representation.

Let $Y$ be an indecomposable $G$-invariant $(Q, {\mathcal
{I}})$-representation, this means that $^{g}Y\cong Y$ for any $g\in
G$. Then, if an indecomposable $(Q, {\mathcal {I}})$-representation
$X$ occurs in summands of $Y$, then all isomorphism classes in
$\{^{g}X\mid g\in G\}$ occur in summands of $Y$ and only once, i.e.,
$Y\cong\bigoplus_{g\in G_{X}}{^{g}X}$.

Finally, for any $G$-invariant $(Q, {\mathcal {I}})$-representation
$Y$, if an indecomposable $(Q, {\mathcal {I}})$-representation $X$
occur in summands of $Y$, then there is a $G$-invariant $(Q, {\mathcal
{I}})$-representation $Y'$ such that $Y\cong Y'\bigoplus
\big(\bigoplus_{g\in G_{X}}{^{g}X}\big)$. Therefore, we can
find finite numbers of indecomposables $(Q, {\mathcal
{I}})$-representations $X_{1}, X_{2}, \cdots,  X_{n}$
such that
$Y\cong \bigoplus_{i=1}^{n} (\oplus_{g\in G_{X_{i}}}{^{g}X_{i}})$
by induction on dimension of $Y$.
The uniqueness follows from the Krull-Schmidt
theorem for finite dimensional $(Q, {\mathcal
{I}})$-representations. \hfill$\square$

\medskip
Now we consider the $\Lambda G$-modules, we have

\medskip
{\bf Lemma 2.2.} Any $\Lambda G$-module $Y$ is a $G$-invariant
$\Lambda$-module.

\medskip
{\bf Proof.} Let $Y$ be a $\Lambda G$-module. We need show that
$^{g}Y\cong Y$ for any $g\in G$. Indeed, for any $g\in G$, we define a map
$f_{g}: {^{g}Y}\rightarrow Y$ by $f_{g}(y)=gy$  for all $~y\in Y$.
Then, for all $a\in \Lambda$ and $y\in Y$, we have
$$f_{g}(a\cdot y)=g(a\cdot y)=g(g^{-1}(a)y)=a(gy)=af_{g}(y).$$
That is, $f_{g}$ is a $\Lambda$-module homomorphism, and hence $f_{g}$ is
an isomorphism. Its inverse is $f_{g^{-1}}: Y\rightarrow {^{g}Y}$ such that
$f_{g^{-1}}(y)=g^{-1}y$ for all $y\in Y$. \hfill$\square$

\medskip
 Now, we consider $\Lambda G$-module $V$ as a
$Q$-representation $(V_{i}, ~V_{a})_{i\in I, a\in A}$.
By Lemma 2.2, we have

(1) If vertices $i$ and $j$ are in the same $G$-orbit, then
$V_{i}\cong V_{j}$ as vector spaces;

(2) If arrows $a$ and $b$ are in the same $G$-orbit, then
 there is a commutative diagram
\[\begin{array}{ccc}
V_{t(a)} & \stackrel{V_{a}}{\longrightarrow} &
V_{h(a)} \\
\llap{$f$}\downarrow  & & \downarrow
\rlap{$g$}\\
V_{t(b)} & \stackrel{V_{b}}{\longrightarrow} & V_{h(b)}
\end{array}\]
such that $f$ and $g$ are isomorphisms. In this case, we denote by
$V_{a}\thickapprox V_{b}$.

\medskip
{\bf Lemma 2.3.} (1) A $\Lambda$-module $Y$ is a $\Lambda G$-module
if and only if for any $g\in G$ there is a $\Lambda$-module homomorphism
$\varphi_{g}: {^{g}Y}\rightarrow Y$ satisfying
$$\varphi_{g}^{|g|}:={^{g^{|g|-1}}\varphi_{g}}\circ\cdots
\circ{^{g}\varphi_{g}}\circ\varphi_{g} =\id_{^{g}Y},$$ where $|g|$
is the order of $g$. In this case, we call that the $\Lambda
G$-module $Y$ is induced by $\{\varphi_{g}\}$.
In particular, each $G$-invariant $\Lambda$-module is a $\Lambda
G$-module.

(2)  For any $g\in G$, if there are two $\Lambda$-module homomorphisms $\varphi_{g},
\psi_{g}: {^{g}Y}\rightarrow Y$, satisfying
$\varphi_{g}^{|g|}=\id_{^{g}Y}$ and $\psi_{g}^{|g|}=\id_{^{g}Y}$
respectively. We denote by $Y_{1}$, $Y_{2}$ the
$\Lambda G$-modules on $Y$ induced by $\{\varphi_{g}\}$ and
$\{\psi_{g}\}$. Then $Y_{1}\cong Y_{2}$ as $\Lambda G$-modules if
and only if there is a $\Lambda$-module isomorphism $f:Y\rightarrow
Y$, such that $f\varphi_{g}=\psi_{g}f$, $\forall ~g\in G$. In
this case, we call $\{\varphi_{g}\}_{g\in G}$ equivalent to
$\{\psi_{g}\}_{g\in G}$.

\medskip
{\bf Proof.} (1)  Assume that $Y$ is a $\Lambda G$-module, we define the
map $\varphi_{g}: {^{g}Y}\rightarrow Y$ by $y\mapsto gy$. Clearly,
$\varphi_{g}^{|g|}=\id_{\,^{g}Y}$, and for any $a\in \Lambda$ and
$y\in Y$, we have
$$\varphi_{g}(ay)=(ga)y=(g(a)g)y=g(a)gy=g(a)\varphi_{g}(y).$$
This means that $\varphi_{g}$ is a $\Lambda$-module homomorphism.
Conversely, if there exists a $\Lambda$-module homomorphism
$\varphi_{g}: {\,^{g}Y}\rightarrow Y$ such that
$\varphi_{g}^{|g|}=\id_{\,^{g}Y}$ for each $g\in G$, we define the
action of $\Lambda G$ on $Y$ by $gay=\varphi_{g}(ay)$ for any $ga\in
\Lambda G$, $y\in Y$. Then $Y$ is a $\Lambda G$-module. Indeed,
$$ga(hby)=\varphi_{g}(a(\varphi_{h}(by)))
=\varphi_{g}(\varphi_{h}(h^{-1}(a)by)) =(ghh^{-1}(a)b)y=(ga\cdot
hb)y,$$ for all $ga, hb \in \Lambda G$, $y\in Y$.

Let $Y$ be a $G$-invariant $\Lambda$-module. That is, there exists a
module isomorphism $\theta_{g}: {^{g}Y}\rightarrow Y$ for every
$g\in G$. As observed in \cite [P.95]{G}, there exists a module
isomorphism $\varphi_{g}: {^{g}Y}\rightarrow Y$ such
that $\varphi_{g}^{|g|}=\id_{^{g}Y}$. So, $Y$ is a $\Lambda
G$-module.

(2) If there is a $\Lambda$-modules isomorphism $f:Y\rightarrow Y$
such that $f\varphi_{g}=\psi_{g}f$ for any $g\in G$. Then
$$f(gay)=f(\varphi_{g}(ay))=\psi_{g}f(ay)=\psi_{g}af(y)=gaf(y),$$
for any $ga\in \Lambda G$ and $y\in Y$, i.e., $f: Y_{1}\rightarrow
Y_{2}$ is a $\Lambda G$-modules isomorphism. Conversely, if there
exists a $\Lambda G$-modules isomorphism $f: Y_{1}\rightarrow
Y_{2}$, one can check that $f$ is a $\Lambda$-modules isomorphism
and satisfy $f\varphi_{g}=\psi_{g}f$ for any $g\in G$.
\hfill$\square$

\medskip
By Lemma 2.2 and 2.3, each $\Lambda$-module $Y$ is a
$\Lambda G$-module if and only if $Y$ is a $G$-invariant
$\Lambda$-module. But for a $G$-invariant $\Lambda$-module $Y$,
It is possible that there are many non-equivalent $\Lambda$-modules
isomorphisms $\{\varphi_{g}\}$ induce $\Lambda G$-module structure
on $Y$.

\medskip
{\bf Question.} How many non-isomorphic $\Lambda G$-module structures
are induced on a $G$-invariant $\Lambda$-module?

\medskip
In the paper \cite{ZL} the answers have been given when $\Lambda=kQ$,
$Q$ is a finite quiver without oriented cycles, and $G$ is a finite
cyclic group. In this section, if $X$ is an indecomposable
$\Lambda$-module satisfying $H_{X}$ is an abelian group,  we will
give  the answer to this question on indecomposable $G$-invariant
$\Lambda$-module $\bigoplus_{g\in G_{X}}{^{g}X}$ (see Theorem 1.1)
by the similar method. Firstly, we need some preparations.

\medskip
The following lemma can be found in \cite{RR} for skew group algebras
of Artin algebras. In fact, it holds for $\Lambda G$ in general.

\medskip
{\bf Lemma 2.4.} Let $X$, $Y$ be indecomposable $\Lambda$-modules
and $G\subseteq \Aut(\Lambda)$, then

(1) $\Lambda G\otimes_{\Lambda}X\cong \bigoplus_{g\in G}{ ^{g}X}$ as
$\Lambda$-modules;

(2) $\Lambda G\otimes_{\Lambda}X\cong \Lambda G\otimes_{\Lambda}Y$
if and only if $Y\cong {^{g}X}$ for some $g\in G$;

(3) The number of summands in the direct sum of indecomposable modules
of $\Lambda G\otimes_{\Lambda}X$ is at most $|H_{X}|$.

\medskip
{\bf Proof.} (1) Note that the subpace
$g\otimes_{\Lambda}X=\{g\otimes x \mid x\in X\}$ of $\Lambda
G\otimes_{\Lambda}X$ has a natural $\Lambda$-module structure given
by $a(g\otimes X)=g\otimes g^{-1}(a)x$.  For any $g\in G$, it is easy to
see that $^{g}X\cong g\otimes_{\Lambda}X$ as $\Lambda$-modules. We have
$\Lambda G\otimes_{\Lambda}X= \bigoplus_{g\in
G}g\otimes_{\Lambda}X\cong \bigoplus_{g\in G}{^{g}X}$.

(2) It is easy to see that $\varphi : \Lambda G\otimes_{\Lambda}X
\rightarrow \Lambda G\otimes_{\Lambda}{^{g}X}$ given by
$\varphi(h\otimes x)=hg\otimes x$ is a $\Lambda G$-module
isomorphism. If $\Lambda G\otimes_{\Lambda}X\cong
\Lambda G\otimes_{\Lambda}Y$, then $Y$ is a summand of $\Lambda
G\otimes_{\Lambda}X\cong \bigoplus_{g\in G}{^{g}X}$ as
$\Lambda$-modules by the statement (1). It follows that
$Y\cong {^{g}X}$ for some $g\in G$ since $Y$ and $^gX$ are indecomposable.

(3) Let $\Lambda G\otimes_{\Lambda}X=Y_{1}\oplus
Y_{2}\oplus\cdots\oplus Y_{t}$, where $Y_{i}$, $i\in
\{1,2,\cdots,t\}$ is an indecomposable $\Lambda G$-module. For each
$i$, we have $Y_{i}\cong\bigoplus_{g\in
S}{^{g}X}$ as $\Lambda$-module, where $S$ is a subset of $G$.
Note that $G_{X}\subseteq S$. We get
 $t\leq |H_{x}|. $\hfill$\square$

\medskip
From now on, we assume that $X$  is an indecomposable $\Lambda$-module
such that $H_{X}=\{g\in G \mid {^{g}X}\cong X\}$ is an abelian group.
It follows that the regular representation $kH_{X}$ can be decomposed
as
$$kH_{X}=\bigoplus_{i=1}^{r}\rho_{i},$$ where $\rho_{i}$ are one dimensional
irreducible $H_{X}$-representations, $r=|H_{X}|$, and
$\rho_{i}\ncong \rho_{j}$ if $i\neq j$.

Since $X$ is an $H_{X}$-invariant $\Lambda$-module, $X$
has a natural $\Lambda H_{X}$-module structure by Lemma 2.3. Thus,
$\rho_{i}\otimes_{k}X$ has also a natural  $\Lambda H_{X}$-module
defined by $ga(l\otimes x)=gl\otimes gax$
for any $ga\in \Lambda H_{X}$ and $l\otimes x\in \rho_{i}\otimes_{k}X$. 
Similarly,  $\Hom_{\Lambda}(X, \rho_{i}\otimes_{k}X)$ is a $\Lambda
H_{X}$-module given by $(gaf)(x)=gaf(x),$  for
$f\in\Hom_{\Lambda}(X, \rho_{i}\otimes_{k}X)$, $ga\in\Lambda H_{X},$
and $x\in X;$ $\rho_{i}\otimes_{k}\End_{\Lambda}(X)$  is $\Lambda
H_{X}$-module given by $ga(l\otimes f)=gl\otimes gaf$ for $h\otimes
f\in \rho_{i}\otimes_{k} \End_{\Lambda}(X)$ and $ga\in \Lambda
H_{X},$ where $(gaf)(x)=gaf(x)$ for each $x\in X$.


\medskip
{\bf Claim.} $\Hom_{\Lambda}(X, \rho_{i}\otimes_{k}X)\cong
\rho_{i}\otimes_{k}\End_{\Lambda}(X)$ as $\Lambda H_{X}$-modules.

Indeed, we define $F_h: \Hom_{\Lambda}(X, \rho_{i}\otimes_{k}X)
\longrightarrow \rho_{i}\otimes_{k}\End_{\Lambda}(X)$ by $f\mapsto
h\otimes \widetilde{f}$,  for each $0\neq h\in \rho_{i}$, where
$\widetilde{f}(x)=ly$ if $f(x)=l\otimes y$, $l\in \rho_{i}$, $x,
y\in X$.  The homomorphism $\widetilde{f}\in \End_{\Lambda}(X)$,
this follows from the fact that $f(ax)=af(x)=a(l\otimes y)=l\otimes
ay$, and hence $\widetilde{f}(ax)=lay=aly=a\widetilde{f}(x)$, for
$a\in \Lambda$, $x\in X$.  This means that $F_{h}$ is well defined.

Now, we show that $F_h$ is a $\Lambda H_{X}$-module homomorphism.
For any $f\in \Hom_{\Lambda}(X, \rho_{i}\otimes_{k}X)$ and $g\in
H_{X}$, we have $(gf)(x)=gf(x)=g(l\otimes y)=gl\otimes gy$. Note
that $gl=\alpha l$ for some $\alpha\in k$ since
$\dim_{k}\rho_{i}=1$. Therefore, $\widetilde{gf}(x)=\alpha
lgy=\alpha gly=\alpha g\widetilde{f}(x)$ and $F_h(gf)=h\otimes
\widetilde{gf}=h\otimes \alpha g \widetilde{f}=\alpha h\otimes g
\widetilde{f}=gh\otimes g
\widetilde{f}=g(h\otimes\widetilde{f})=gF_h(f)$. So $F_h$ is a
$kH_{X}$-module homomorphism. Moreover, since
$(af)(x)=af(x)=a(l\otimes y)=l\otimes ay$, $\widetilde{af}(x)=lay
=aly=a\widetilde{f}(x)$ for all $f\in \Hom_{\Lambda}(X,
\rho_{i}\otimes_{k}X)$, $a\in \Lambda$, $F_h$ is indeed a
$\Lambda$-module homomorphism.

Finally, noting that $F$ is injective and $\dim_{k}\Hom_{\Lambda}(X,
\rho_{i}\otimes_{k}X)=\dim_{k}\rho_{i}\otimes_{k}\End_{\Lambda}(X)$,
 we have $F$ is a $\Lambda H_{X}$-module isomorphism.

For convenience, let us denote

{\bf Condition C}: $X$ is an indecomposable $\Lambda$-module and
$H_{X}=\{g\in G \mid {^{g}X}\cong X\}$ is abelian.
$kH_{X}=\bigoplus_{i=1}^{r}\rho_{i}$ is the irreducible
decomposition of $kH_{X}$ as $H_{X}$-representations, $r=|H_{X}|$.

\medskip
{\bf Lemma 2.5.}  Assume that  {\bf Condition C}. Then for any $i\in
\{1,2,\cdots,r\}$, we have

(1) $\rho_{i}\otimes_{k}X\cong X$ as $\Lambda$-modules and
$\rho_{i}\otimes_{k}X$ is indecomposable as a $\Lambda H_{X}$-module
for  all $~i\in \{1,2,\cdots,r\}$;

(2) $\rho_{i}\otimes_{k}X \ncong \rho_{j}\otimes_{k}X$ as $\Lambda
H_{X}$-modules, if $i\neq j$.

(3) $\Lambda H_{X}\otimes _{\Lambda}X\cong
\bigoplus_{i=1}^{r}\rho_{i}\otimes_{k}X$ as $\Lambda H_{X}$-modules;

(4) For any $\Lambda H_{X}$-module $Y$, if $Y\cong X$ as
$\Lambda$-modules, then there exists a unique $i\in
\{1,2,\cdots,r\}$ such that $Y\cong \rho_{i}\otimes_{k}X$ as
$\Lambda H_{X}$-modules. Hence there are $r$ non-isomorphic $\Lambda
H_{X}$-modules induced from $X$.

\medskip
{\bf Proof.} (1) For each $0\neq l\in \rho_{i}$, we define a
bijection $f : X \rightarrow \rho_{j}\otimes_{k}X$ by $x\mapsto
l\otimes x$. Observe that $f$ is a $\Lambda$-module homomorphism
since $f(ax)=l\otimes ax=a(l\otimes x)=af(x)$ for all $a\in
\Lambda$, $x\in X$. Therefore, $\rho_{i}\otimes_{k}X$ is an
indecomposable $\Lambda$-module, and hence an indecomposable
$\Lambda H_{X}$-module.

(2) Assume that $\rho_{i}\otimes_{k}X \cong \rho_{j}\otimes_{k}X$.
We have $\rho_{i}\otimes_{k}\End_{\Lambda}(X)\cong
\rho_{j}\otimes_{k}\End_{\Lambda}(X)$ by the claim. Since
$\End_{\Lambda}(X)/\rad\End_{\Lambda}(X)\cong k$ and
$\rad\End_{\Lambda}(X)$ is closed under the action of $H_{X}$, we
have
$$\rho_{i}\otimes_{k}\End_{\Lambda}(X)/\rad\End_{\Lambda}(X)\cong
\rho_{j}\otimes_{k}\End_{\Lambda}(X)/\rad\End_{\Lambda}(X).$$ This
means $\rho_{i}\cong \rho_{j}$ as $H_{X}$-modules and we get a
contradiction.

(3) By \cite[Lemma 3.2.1]{L}, $\rho_{i}\otimes_{k}X\mid\Lambda
H_{X}\otimes _{\Lambda}(\rho_{i}\otimes_{k}X)$. This induce that
$\rho_{i}\otimes_{k}X\mid\Lambda H_{X}\otimes _{\Lambda}X$ since
$\rho_{i}\otimes_{k}X\cong X$ as $\Lambda$-modules. Note that
$\rho_{i}\otimes_{k}X\ncong \rho_{j}\otimes_{k}X$ if $i\neq j$, we
have $\big(\bigoplus_{i=1}^{r}\rho_{i}\otimes_{k}X\big)\mid\Lambda
H_{X}\otimes _{\Lambda}X$. That is, $\Lambda H_{X}\otimes
_{\Lambda}X\cong \bigoplus_{i=1}^{r}\rho_{i}\otimes_{k}X$ by Lemma
2.4(3).

(4) If $Y$ is a $\Lambda H_{X}$-module, and $Y\cong X$ as
$\Lambda$-modules, then $Y$ is an indecomposable $\Lambda
H_{X}$-module. Since $Y\mid \Lambda
H_{X}\otimes_{\Lambda}Y\cong\Lambda H_{X}\otimes_{\Lambda}X$, it is easy to
see that there exists a unique $i\in \{1,2,\cdots,r\}$ such that $Y\cong
\rho_{i}\otimes_{k}X$. \hfill$\square$

\medskip
{\bf Lemma 2.6.}  Assume  that {\bf Condition C}. Then we have
for any $i\in \{1, 2, \cdots, r\}$

(1) $\Lambda G\otimes_{\Lambda
H_{X}}(\rho_{i}\otimes_{k}X)\cong\bigoplus_{g\in G_{X}}{^{g}X}$ as
$\Lambda$-modules;

(2) $\Lambda G\otimes_{\Lambda H_{X}}(\rho_{i}\otimes_{k}X)$ is an
indecomposable $\Lambda G$-module;

(3) $\Lambda G\otimes_{\Lambda
H_{X}}(\rho_{i}\otimes_{k}X)\ncong\Lambda G\otimes_{\Lambda
H_{X}}(\rho_{j}\otimes_{k}X)$ as $\Lambda G$-modules, if $i\neq j$;

(4) $\Lambda G\otimes_{\Lambda}X\cong\bigoplus_{i=1}^{r}\Lambda
G\otimes_{\Lambda H_{X}}(\rho_{i}\otimes_{k}X)$ as $\Lambda
G$-modules;

\medskip
{\bf Proof.} (1) Since $\Lambda G\otimes_{\Lambda
H_{X}}(\rho_{i}\otimes_{k}X)\cong\bigoplus_{g\in G_{X}}g\otimes
\rho_{i}\otimes X$ as $\Lambda H_{X}$-modules and $\rho_{i}\otimes
X\cong X$ as $\Lambda$-modules, we have $\Lambda G\otimes_{\Lambda
H_{X}}(\rho_{i}\otimes_{k}X)\cong\bigoplus_{g\in G_{X}}g\otimes
X\cong\bigoplus_{g\in G_{X}}{^{g}X}$.

(2) The result follows  from that $\Lambda G\otimes_{\Lambda
H_{X}}(\rho_{i}\otimes_{k}X)\cong\bigoplus_{g\in G_{X}}{^{g}X}$ is
an
indecomposable $G$-invariant $\Lambda$-module. 

(3) Suppose that $\Lambda G\otimes_{\Lambda
H_{X}}(\rho_{i}\otimes_{k}X)\cong\Lambda G\otimes_{\Lambda
H_{X}}(\rho_{j}\otimes_{k}X)$. Note that
$$\Lambda G\otimes_{\Lambda H_{X}}(\rho_{i}\otimes_{k}X)\cong
\bigoplus_{g\in G_{X}}g\otimes \rho_{i}\otimes X$$ as $\Lambda
H_{X}$-modules, we have $h\otimes \rho_{i}\otimes X\mid\Lambda
G\otimes_{\Lambda H_{X}}(\rho_{j}\otimes_{k}X)\cong\bigoplus_{g\in
G_{X}}g\otimes \rho_{j}\otimes X$. If $h\otimes \rho_{i}\otimes
X\cong h\otimes \rho_{j}\otimes X$, then $\rho_{i}\otimes X\cong
\rho_{j}\otimes X$ as $\Lambda H_{X}$-modules. This contradict to
Lemma 2.5. If $h\otimes \rho_{i}\otimes X\cong g\otimes
\rho_{j}\otimes X$ for some $h\ne g\in G_{X}$, we have
$X\cong{^{g}X}$ as $\Lambda$-modules since $g\otimes \rho_{j}\otimes
X\cong {^{g}X}$. This is also a contradiction.

(4)  Note that $\Lambda G\otimes_{\Lambda
H_{X}}(\rho_{i}\otimes_{k}X)\mid\Lambda G\otimes_{\Lambda}\Lambda
G\otimes_{\Lambda H_{X}}(\rho_{i}\otimes_{k}X)$, by the statement
(1) and Lemma 2.4(2) we have $\Lambda G\otimes_{\Lambda
H_{X}}(\rho_{i}\otimes_{k}X)\mid\Lambda
G\otimes_{\Lambda}\big(\bigoplus_{g\in G_{X}}{^{g}X}\big)$ and
$\Lambda G\otimes_{\Lambda H_{X}}(\rho_{i}\otimes_{k}X)\mid\Lambda
G\otimes_{\Lambda}X$ for any $i\in \{1,2,\cdots,r\}$. Thus,
$\big(\bigoplus_{i=1}^{r}\Lambda G\otimes_{\Lambda
H_{X}}(\rho_{i}\otimes_{k}X)\big)\mid\Lambda G\otimes_{\Lambda}X$,
and $\Lambda G\otimes_{\Lambda}X\cong\bigoplus_{i=1}^{r}\Lambda
G\otimes_{\Lambda H_{X}}(\rho_{i}\otimes_{k}X)$ by Lemma 2.4(3).
\hfill$\square$

\medskip
{\bf Proof of Theorem 1.1.} By Lemma 2.1 and 2.2, if $\Lambda
G$-module $Y$ satisfies $Y\cong\bigoplus_{g\in G_{X}}{^{g}X}$ as
$\Lambda$-modules, then $Y$ is an indecomposable $\Lambda G$-module.
By \cite [Lemma 3.2.1]{L}, $Y\mid \Lambda G\otimes_{\Lambda}Y\cong
\Lambda G\otimes_{\Lambda}\big(\bigoplus_{g\in G_{X}}{^{g}X}\big)$.
On the other hand,
 $\Lambda G\otimes_{\Lambda}{^{g}X}\cong\Lambda
G\otimes_{\Lambda}{X}$ for any $g\in G$, we have $Y\mid \Lambda
G\otimes_{\Lambda}X$. Thus there exists a unique $i\in
\{1,2,\cdots,r\}$ such that $Y\cong \Lambda G\otimes_{\Lambda
H_{X}}(\rho_{i}\otimes_{k}X)$, by Lemma 2.6(4) and Krull-Schmidt
Theorem. \hfill$\square$

\medskip
{\bf Proof of Theorem 1.2.} Suppose $Y$ is an indecomposable
$\Lambda G$-module. Then $Y$ is a $G$-invariant $\Lambda$-module and
$Y\cong\bigoplus_{j=1}^{s}\big(\bigoplus_{g\in
G_{X_{j}}}{^{g}X_{j}}\big)$ with some indecomposable
$\Lambda$-modules $X_{1}, X_{2}, \cdots, X_{s}$. Since $Y\mid\Lambda
G\otimes_{\Lambda}Y\cong \bigoplus_{j=1}^{s}\bigoplus_{g\in
G_{X_{j}}}\Lambda G\otimes_{\Lambda}{^{g}X_{j}}$, there exists a$j$
such that $Y\mid\Lambda G\otimes_{\Lambda}{X_{j}}$. If we denote by
$kH_{X_{j}}=\bigoplus_{i=1}^{r_{j}}\rho_{i}^{j}$ the irreducible
decomposition of $kH_{X_{j}}$ as $H_{X_{j}}$-representations, then
there exists a unique $\rho_{i}^{j}$ such that $Y\cong \Lambda
G\otimes_{\Lambda H_{X_{j}}}(\rho_{i}^{j}\otimes_{k}X_{j})$ and
$Y\cong\bigoplus_{g\in G_{X_{j}}}{^{g}X_{j}}$ as $\Lambda$-modules,
by Lemma 2.6. Thus, $Y$ is an indecomposable $G$-invariant
$\Lambda$-module.

If $Y$ is an indecomposable $G$-invariant $\Lambda$-module. By Lemma
2.3, $Y$ is a $\Lambda G$-module and  indecomposable since
any $\Lambda G$-module is a $G$-invariant $\Lambda$-module.
\hfill$\square$

\medskip
According to Theorem 1.1 and Theorem 1.2, we can obtain all $\Lambda
G$-modules from indecomposable $\Lambda$-modules under the
assumption that $G$ is abelian. In this case, for any indecomposable
$\Lambda$-module $X$, the $G$-invariant $\Lambda$-module
$\bigoplus_{g\in G_{X}}{^{g}X}$ has $r$ non-isomorphic $\Lambda
G$-module structures, where $r=|H_X|$. It is equivalent to say that
there are $r$ non-equivalent $\Lambda$-module isomorphisms
$\{\varphi_{g}: {^{g}X}\rightarrow X\}$. Moreover,  we have

\medskip
{\bf Proposition 2.7.} Suppose $G$ is an abelian group. For any
indecomposable $\Lambda G$-module $Y$, we have

(1) $Y$ is simple if and only if there exists a simple
$\Lambda$-module $\mathscr{S}$, such that $Y\cong \bigoplus_{g\in
G_{\mathscr{S}}}{^{g}{\mathscr{S}}}$.

(2) $Y$ is projective if and only if there exists an indecomposable
projective $\Lambda$-module $\mathscr{P}$, such that $Y\cong
\bigoplus_{g\in G_{\mathscr{P}}}{^{g}{\mathscr{P}}}$.

(2) $Y$ is injective if and only if there exists an indecomposable
injective $\Lambda$-module $\mathscr{I}$, such that $Y\cong
\bigoplus_{g\in G_{\mathscr{I}}}{^{g}{\mathscr{I}}}$.

\medskip
{\bf Proof.} Suppose $Y=\bigoplus_{g\in
G_{{\mathscr{S}}}}{^{g}{\mathscr{S}}}$ for some simple
$\Lambda$-module $\mathscr{S}$. If there is a proper submodule $Y'$ of
$Y$, then $Y'$ is a summand of $Y$ as a $\Lambda$-module and
$Y'=\bigoplus_{g\in E}{^{g}{\mathscr{S}}}$, where $E$ is a proper
subset of $G_{{\mathscr{S}}}$. This is a  contradiction  by Lemma 2.1
since $Y'$ is $G$-invariant. Conversely, Given a simple $\Lambda
G$-module $Y$, by Theorem 1.2, there is an indecomposable
$\Lambda$-module $X$, such that $Y\cong\bigoplus_{g\in
G_{X}}{^{g}X}$. If $X$ is not a simple $\Lambda$-module, there exists a
proper submodule $X'$ of $X$. It follows that
$\bigoplus_{g\in G_{X'}}{^{g}X'}$ is a proper submodule of $Y$. This is also
a contradiction and the statement (1) is proved.

Since $X$ is a projective $\Lambda G$-module if and only if $X$ as
$\Lambda$-module is projective (cf. \cite [Lemma 3.1.7]{L}), we get
(2). By duality, we have (3). \hfill$\square$

\medskip
Finally, as an example, we consider the skew group algebras of
deformed preprojective algebras. For each finite group $G\subseteq
\Aut(Q)$, $G$ acts in a natural way on the double quiver
$\overline{Q}$, i.e., $g(a')=g(a)'$ for all $g\in G$ and additive
arrow $a'$.

The deformed preprojective algebras defined by
$\Pi_{Q}^{\lambda}:=k\overline{Q}/{\mathcal {I}}$, where
$\lambda=(\lambda_{i})_{i\in I}\in k^{I}$, $\mathcal {I}$ is the
ideal of path algebra $k\overline{Q}$ generated by the deformed
preprojective relations. Note that $g({\mathcal {I}})={\mathcal
{I}}$ for any $g\in G$ if and only if $\lambda=(\lambda_{i})_{i\in
I}$ satisfies
$$(\ast) : \qquad \lambda_{i}=\lambda_{j} \mbox{ if vertices } i \mbox{
and } j \mbox{ are in the same }
G\mbox{-orbit}.\qquad\qquad\qquad\qquad\qquad\qquad$$ Thus, if
$\lambda=(\lambda_{i})_{i\in I}$ satisfies ($\ast$), then
$G\subseteq \Aut(\overline{Q}, {\mathcal {I}})$ and we  get

\medskip
{\bf Corollary 2.8.} If $G\subseteq \Aut(\overline{Q}, {\mathcal
{I}})$ is abelian and $\lambda=(\lambda_{i})_{i\in I}$ satisfies
($\ast$), then a $\Pi_{Q}^{\lambda}$-module $Y$ is an
indecomposable (projective, injective) $\Pi_{Q}^{\lambda} G$-module
if and only if $Y$ is an indecomposable (projective, injective)
$G$-invariant $\Pi_{Q}$-module.
Moreover, for any $G$-invariant $\Pi_{Q}$-module $Y$, the number of
non-isomorphic $\Pi_{Q}^{\lambda} G$-module structure on $Y$ can be
determined. \hfill$\square$

\section{The reflection functors}
In this section, we introduce a reflection functor for the
module category of skew group algebras of deformed preprojective
algebras.

\medskip
Firstly, let us recall some notations. The dimension vector of any
representation $V=(V_{i}, ~V_{a})_{i\in I, a\in A}$ of quiver $Q$ is
denoted by ${\bf dim}\,V=(\dim_{k}V_{i})_{i\in I}$. The Ringel form for $Q$ is
defined to be the bilinear form on ${\mathbb{Z}}^{I}$ with
$$\langle\alpha, \beta\rangle=\sum_{i\in I}\alpha_{i}\beta_{i}-
\sum_{a\in A}\alpha_{t(a)}\beta_{h(a)},$$ for $\alpha, \beta \in
{\mathbb{Z}}^{I}$. The bilinear form $(\alpha,
\beta)=\langle\alpha, \beta\rangle+\langle\beta, \alpha\rangle$
is the corresponding symmetric bilinear form.

We say that a vertex $i$ is loop-free if there are no arrows
$a:i\rightarrow i$, and if so, we define simple reflection
$r_{i}:{\mathbb{Z}}^{I}\rightarrow {\mathbb{Z}}^{I}$ by
$$r_{i}(\alpha)=\alpha-(\alpha, \varepsilon_{i})\varepsilon_{i},$$
where $\varepsilon_{i}$ is the coordinate vector at $i$. There is a
dual reflection $d_{i}:k^{I}\rightarrow k^{I}$ define by
$(d_{i}\lambda)_{j}=\lambda_{j}-(\varepsilon_{i},
\varepsilon_{j})\lambda_{i}$, for any $j\in I$, $\lambda \in k^{I}$.
By direct calculation, we have $d_{i}\lambda\cdot\alpha=\lambda\cdot
r_{i}\alpha$ for any $\lambda \in k^{I}$ and $\alpha\in
{\mathbb{Z}}^{I}$.

Let $G\subseteq \Aut(Q)$  be a finite group. For $a\in A$, $i\in I$, we
denote by ${\mathcal {O}}_{a}$ and ${\mathcal {O}}_{i}$ the orbit of
$a$ and $i$ respectively under the action of $G$, by $S_{{\mathcal {O}}_{a}}Q$
the quiver obtained from $Q$ by reversing the direction of all
arrows in ${\mathcal {O}}_{a}$, and by $S_{{\mathcal {O}}_{i}}Q$ the
quiver obtained from $Q$ by reversing the direction of all arrows
$a$ satisfying $h(a)$ or $t(a)$ in ${\mathcal {O}}_{i}$,
respectively.
Obviously, for any $i\in I$, there exist $\{a_{1}, a_{2}, \cdots,
a_{n}\}\subseteq A$ such that ${\mathcal
{O}}_{i}=\bigcup_{j=1}^{n}{\mathcal {O}}_{a_{j}}$ and ${\mathcal
{O}}_{a_{j}}\neq {\mathcal {O}}_{a_{l}}$ if $j\neq l$. Note that for
any arrows $a, a'\in A$, $a'$ is also a arrow in $S_{{\mathcal
{O}}_{a}}Q$ if $a'\notin{\mathcal
{O}}_{a}$, we have $$S_{{\mathcal {O}}_{i}}Q=S_{{\mathcal
{O}}_{a_{n}}}\cdots S_{{\mathcal {O}}_{a_{2}}}S_{{\mathcal
{O}}_{a_{1}}}Q.$$
If we denote by $b^{\ast}$ the reversed arrow in
$S_{{\mathcal {O}}_{a}}Q$ corresponding to $b$, the
action of $G$ on $S_{{\mathcal {O}}_{a}}Q$ is defined by
$g(b^{\ast})=g(b)^{\ast}$, then $G\subseteq \Aut(S_{{\mathcal
{O}}_{a}}Q)$ and hence $G\subseteq \Aut(S_{{\mathcal {O}}_{i}}Q)$, for all
$a\in A$ and $i\in I$.

The action of $G$ on $Q$ is called admissible if $Q$ have
no arrows connecting two vertices in the same $G$-orbit.
If the action of $G$  is admissible, then
$$r_{{\mathcal {O}}_{i}}=\prod_{j\in {\mathcal {O}}_{i}}r_{j}
\mbox{~~ and ~~} d_{{\mathcal {O}}_{i}}=\prod_{j\in {\mathcal
{O}}_{i}}d_{j}$$ are well defined.

\medskip
Considering the skew group algebra $\Pi_{Q}^{\lambda} G$, where
$\lambda=(\lambda_{i})_{i\in I}$ satisfying $(\ast):$
$\lambda_{i}=\lambda_{j}$ if vertices $i$ and $j$ are in the same
orbit under the action of $G$. If the action of $G$ is admissible,
$d_{{\mathcal{O}}_{i}}\lambda$ also satisfies $(\ast)$ . Therefore,
we get a new skew group algebra $\Pi_{Q}^{d_{{\mathcal
{O}}_{i}}\lambda} G$. By \cite [Lemma 2.2]{CH}, it is well known
that a deformed preprojective algebra $\Pi_{Q}^{\lambda}$ does not
depend on the orientation of $Q$. Therefore, we have algebraic
isomorphisms $\Pi_{Q}^{\lambda}\cong \Pi_{S_{{\mathcal
{O}}_{a}}Q}^{\lambda}\cong \Pi_{S_{{\mathcal {O}}_{i}}Q}^{\lambda}$
for any $a\in A$, $i\in I$.

\medskip
{\bf Lemma 3.1.}
$\Pi_{Q}^{\lambda}G\cong \Pi_{S_{{\mathcal
{O}}_{a}}Q}^{\lambda}G\cong \Pi_{S_{{\mathcal
{O}}_{i}}Q}^{\lambda}G$ as algebras for any $a\in A$, $i\in I$.

\medskip
{\bf Proof.} Note that the action of $G$ on $S_{{\mathcal
{O}}_{a}}Q$ is given by $g(b^{\ast})=g(b)^{\ast}$, for any $g\in G$
and any reversed arrow $b^{\ast}$ of $b$ in $S_{{\mathcal {O}}_{a}}Q$, there
is an isomorphism $F:\Pi_{Q}^{\lambda}G\rightarrow \Pi_{S_{{\mathcal
{O}}_{a}}Q}^{\lambda}G$ by  sending $b$ to $(b^{\ast})'$ and $b'$ to
$-b^{\ast}$, for all $b\in {\mathcal {O}}_{a}$.  We conclude that
$\Pi_{Q}^{\lambda}G\cong \Pi_{S_{{\mathcal {O}}_{a}}Q}^{\lambda}G$,
and $\Pi_{Q}^{\lambda}G\cong \Pi_{S_{{\mathcal
{O}}_{i}}Q}^{\lambda}G$. \hfill$\square$

\medskip
In \cite{CH}, Crawley-Boevey and Holland introduced an interesting
reflection functor for deformed preprojective algebras. We want to
define a reflection functor for skew group
algebras of deformed preprojective algebras.

Firstly, we recall the definition of the reflection functor defined
in \cite{CH}. Let $V$ be a $\Pi_{Q}^{\lambda}$-module. We can
identify $V$ with a representation $(V_{i}, ~V_{a})_{i\in I, a\in
\overline{A}}$ of double quiver $\overline{Q}$ satisfying the
deformed preprojective relations
$$\sum_{a\in A, ~h(a)=j}V_{a}V_{a'}- \sum_{a\in A,
~t(a)=j}V_{a'}V_{a}-\lambda_{j}\id_{V_{j}}, \qquad \forall j\in I.$$
Suppose that $i\in I$ is loop-free, $\lambda_{i}\neq 0$, and no
$a\in A$ such that $t(a)=i$. We define a new representation
$E_{i}V=W:=(W_{i}, ~W_{a})_{i\in I, a\in \overline{A}}$ for the double
quiver $\overline{Q}$ by let $W_{i}=\Ker \pi$, $W_{i}=V_{j}$ for
$j\neq i$, and with linear maps $W_{a}=V_{a}$, $W_{a'}=V_{a'}$ if
$h(a)\notin {\mathcal {O}}_{i}$, while if $h(a)=j\in {\mathcal
{O}}_{i}$,
$$W_{a}=-\lambda_{j}(1-\mu\pi)\mu_{a} :
W_{t(a)}\rightarrow W_{j}, ~~~~ W_{a'}=\pi_{a}|_{W_{j}} :
W_{j}\rightarrow W_{t(a)},$$ where $V_{\oplus}=\bigoplus_{a\in A,
~h(a)=j}V_{t(a)}$, $\mu_{a} : V_{t(a)}\rightarrow V_{\oplus}$,
~$\pi_{a} : V_{\oplus}\rightarrow V_{t(a)}$ is the canonical
inclusion and projection, and $\mu=\sum_{a\in A,
~h(a)=j}\mu_{a}V_{a'}$, $\pi=\frac{1}{\lambda_{j}}\sum_{a\in A,
~h(a)=j}V_{a}\pi_{a}$. One can check that $W=(W_{i}, ~W_{a})_{i\in
I, a\in \overline{A}}$ is a $\Pi_{Q}^{d_{i}\lambda}$-module, and
$E_{i}V=E_{i}V'\oplus E_{i}V''$ if $V=V'\oplus V''$. Therefore,
$E_{i}V$ is indecomposable if and only if $V$ is indecomposable.

\medskip
{\bf Theorem 3.2 $\cite{CH}$.} If $i$ is a loop-free vertex,
$\lambda_{i}\neq 0$, then there is an equivalence
$$\begin{aligned}E_{i}:\qquad\Pi_{Q}^{\lambda}\mbox{-{\bf mod}}\qquad
&\longrightarrow
\qquad\Pi_{Q}^{d_{i}\lambda}\mbox{-{\bf mod}} \nonumber\\
 V=(V_{i}, ~V_{a})_{i\in I, a\in
\overline{A}}  &\mapsto  ~~W=(W_{i}, ~W_{a})_{i\in I, a\in
\overline{A}}
\end{aligned}$$ that acts as the
reflection $r_{i}$ on dimension vector. \hfill$\square$

\medskip
Now, we consider the action of $G$. If the action of $G$ on $Q$ is
admissible and $\lambda_{i}\neq 0$, then functor $E_{{\mathcal
{O}}_{i}}:=\prod_{j\in {\mathcal {O}}_{i}}E_{j} :
\Pi_{Q}^{\lambda}\mbox{-{\bf mod}}\rightarrow \Pi_{Q}^{d_{{\mathcal
{O}}_{i}}\lambda}\mbox{-{\bf mod}}$ is well defined, since there is
no arrows in $Q$ connecting two vertices in the same $G$-orbit and
$\lambda_{j}\neq 0$ for each $j\in {\mathcal {O}}_{i}$. Moreover, it
is easy to see that $E_{{\mathcal {O}}_{i}}$ is  a reflection
functor since so $E_{j}$, $j\in {\mathcal {O}}_{i}$ are.

\medskip
{\bf Lemma 3.3.} Assume that $i$ is a loop-free vertex and $\lambda_{i}\neq 0$,
and the action of $G$ is admissible, then for each indecomposable
$\Pi_{Q}^{\lambda}$-module $X$, we have

(1) $H_{X}=H_{E_{{\mathcal {O}}_{i}}X}$. In this case, we can take
$G_{E_{{\mathcal {O}}_{i}}X}=G_{X}$.

(2)  $E_{{\mathcal {O}}_{i}}{^{g}X}\cong {^{g}(E_{{\mathcal
{O}}_{i}}X)}$ as $\Pi_{Q}^{d_{{\mathcal {O}}_{i}}\lambda}$-modules,
for each $g\in G$. Thus
$$E_{{\mathcal {O}}_{i}}\big(\bigoplus_{g\in
G_{X}}{^{g}X}\big)\cong \bigoplus_{g\in G_{E_{{\mathcal
{O}}_{i}}X}}{^{g}(E_{{\mathcal {O}}_{i}}X)}.$$

\medskip
{\bf Proof.} We can identify $X$ with representations $(X_{i}, ~X_{a})_{i\in I, a\in
\overline{A}}$, and $E_{{\mathcal
{O}}_{i}}X$ with $(Y_{i}, ~Y_{a})_{i\in I, a\in \overline{A}}$, which are
representations of double quiver $\overline{Q}$ satisfying the corresponding
preprojective relations, respectively.

(1) For each $g\in H_{X}$, we have $X\cong{ ^{g}X}$. This means that
$X_{i}\cong X_{g^{-1}(i)}$, $X_{a}\approx X_{g^{-1}(a)}$ as well as
$X_{a'}\approx X_{g^{-1}(a')}$ for any $i\in I$, $a\in A$. By the
definition of $E_{i}$, it is easy to see that $Y_{i}\cong
Y_{g^{-1}(i)}$, $Y_{a}\approx Y_{g^{-1}(a)}$ and $Y_{a'}\approx
Y_{g^{-1}(a')}$ for any $i\in I$, $a\in A$. Hence $g\in
H_{E_{{\mathcal {O}}_{i}}X}$ and $H_{X}\subseteq H_{E_{{\mathcal
{O}}_{i}}X}$. Similarly, we have
$H_{E_{{\mathcal{O}}_{i}}X}\subseteq H_{X}$ since
$E_{{\mathcal {O}}_{i}}$ is a reflection functor.
Therefore $H_{X}= H_{E_{{\mathcal {O}}_{i}}X}$.

(2) Note that  $E_{j}{^{g}X}\cong { {^{g}(E_{g(j)}X)}}$ as
$\Pi_{Q}^{d_{j}\lambda}$-modules for all $j\in {\mathcal {O}}_{i}$,
$g\in G$, and $g({\mathcal {O}}_{i})={\mathcal {O}}_{i}$. Accordingly,
if  ${\mathcal {O}}_{i}:=\{j_{1}, j_{2},\cdots, j_{n}\}$,
it is easy to
see that
$$E_{{\mathcal {O}}_{i}}{^{g}X}=E_{j_{n}}\cdots
E_{j_{2}}E_{j_{1}}{^{g}X}\cong {^{g}(E_{g(j_{n})}\cdots
E_{g(j_{2})}E_{g(j_{1})}X)}={^{g}(E_{{\mathcal {O}}_{i}}X)}.$$
Since the reflection functor $E_{i}$ preserve direct sum, it follows
that
$$E_{{\mathcal {O}}_{i}}\big(\bigoplus_{g\in
G_{X}}{^{g}X}\big)\cong \bigoplus_{g\in G_{X}}{^{g}(E_{{\mathcal
{O}}_{i}}X)}=\bigoplus_{g\in G_{E_{{\mathcal
{O}}_{i}}X}}{^{g}(E_{{\mathcal {O}}_{i}}X)}.$$
The proof of the lemma is compeleted.\hfill$\square$

\medskip
{\bf Proof of Theorem 1.3.} It is sufficient to prove the
theorem  for indecomposable modules.

By Theorem
1.1 and Theorem 1.2,
 there is an
indecomposable $\Pi_{Q}^{\lambda}$-module $X$ such that $V\cong
\bigoplus_{g\in G_{X}}{^{g}X}$  for any indecomposable
$\Pi_{Q}^{\lambda}G$-module $V$. Note that the
$\Pi_{Q}^{\lambda}G$-module structure on $\bigoplus_{g\in
G_{X}}{^{g}X}$ is induced by some $\{\varphi_{g}\}$, where
$\varphi_{g}: {X^{g}}\rightarrow X$ is $\Pi_{Q}^{\lambda}$-modules
isomorphism and satisfies $\varphi_{g}^{|g|}=\id_{^{g}Y}$ for any
$g\in G$. In this case, the $\Pi_{Q}^{\lambda}G$-module $V$ is said
to be induced by $\{\varphi_{g}\}$.

Given an indecomposable $\Pi_{Q}^{\lambda}G$-module $V=
\bigoplus_{g\in G_{X}}{^{g}X}$ induced by $\{\varphi_{g}\}$.
Note that $E_{{\mathcal{O}}_{i}}$ is the reflection functor
of $\Pi_{Q}^{\lambda}$-modules,
we have $E_{{\mathcal {O}}_{i}}\varphi_{g}:{^{g}(E_{{\mathcal
{O}}_{i}}X)}\rightarrow E_{{\mathcal {O}}_{i}}X$ is an
 isomorphism as $\Pi_{Q}^{\lambda}$-modules and $(E_{{\mathcal
{O}}_{i}}\varphi_{g})^{|g|}=\id_{{^{g}(E_{{\mathcal {O}}_{i}}X)}}$
for any $g\in G$. Now, we define
$F_{i}V=\bigoplus_{g\in G_{X}}{^{g}(E_{{\mathcal {O}}_{i}}X)}.$
Note that the $\Pi_{Q}^{d_{{\mathcal {O}}_{i}}\lambda}G$-module
$F_{i}V$ is induced by $E_{{\mathcal{O}}_{i}}\varphi_{g}$.

Let us suppose that the $\Pi_{Q}^{\lambda}G$-modules $V$ and $V'$
are induced by $\{\varphi_{g}\}$ and $\{\psi_{g}\}$ respectively.
For any homomorphism $f : V\rightarrow V'$ in
$\Pi_{Q}^{\lambda}G\mbox{-{\bf mod}}$, we define the homomorphism
$F_{i}f:=E_{{\mathcal {O}}_{i}}f : F_{i}V\rightarrow F_{i}V'$. The
functor $F_{i}$ is well defined. Indeed

$$\begin{aligned}
F_{i}f(gav)=&E_{{\mathcal
{O}}_{i}}f(E_{{\mathcal {O}}_{i}}\varphi_{g}(av))=E_{{\mathcal
{O}}_{i}}(f\varphi_{g})(av)=E_{{\mathcal {O}}_{i}}(\psi_{g}f)(av)
\nonumber\\
=&E_{{\mathcal {O}}_{i}}\psi_{g}(E_{{\mathcal {O}}_{i}}f(av))
=E_{{\mathcal {O}}_{i}}\psi_{g}(aE_{{\mathcal {O}}_{i}}f(v))
=gaF_{i}f(v),\end{aligned}$$ for any $ga\in \Pi_{Q}^{\lambda}G$,
$v\in F_{i}V$. Thus, we obtain
the functor
$$F_{i}:\Pi_{Q}^{\lambda}G\mbox{-{\bf mod}}\rightarrow
\Pi_{Q}^{d_{{\mathcal {O}}_{i}}\lambda}G\mbox{-{\bf mod}}.$$
Moreover, since $F_{i}V\cong E_{{\mathcal {O}}_{i}}V$ as
$\Pi_{Q}^{\lambda}$-modules, we have ${\bf dim}F_{i}V=r_{{\mathcal
{O}}_{i}}{\bf dim}V$ by Theorem 3.2.

As a matter of fact, we can also define a functor
$$F'_{i}:\Pi_{Q}^{d_{{\mathcal {O}}_{i}}\lambda}G\mbox{-{\bf
mod}}\longrightarrow \Pi_{Q}^{\lambda}G\mbox{-{\bf mod}}$$ in
a similar way. Note that $H_{X}=H_{E_{{\mathcal {O}}_{i}}X}$,
it is straightforward to check that there are natural isomorphisms
$V\to F'_{i}F_{i}V$
and $W\to F_{i}F'_{i}W$. Therefore $F_{i}$ is an equivalence.
\hfill$\square$

\medskip
Now, we denote by $\mathcal {W}$ the group generated by the
reflections $r_{{\mathcal {O}}_{i}}$ for all loop-free vertices
$i\in I$. The length $l(w)$ of $w\in {\mathcal {W}}$ is said to be $n$
 if
there is a expression $w=r_{{\mathcal {O}}_{i_{n}}}\cdots
r_{{\mathcal {O}}_{i_{2}}}r_{{\mathcal {O}}_{i_{1}}}$ with shortest
length. For each $w\in
{\mathcal {W}}$, we define an action of $w$ on $k^{I}$ by
$$(w\lambda)\cdot\alpha=\lambda\cdot(w^{-1}\alpha)$$
for all $\lambda\in k^{I}$, $\alpha\in {\mathbb{Z}}^{I}$. If
$w=r_{{\mathcal {O}}_{i_{n}}}\cdots r_{{\mathcal
{O}}_{i_{2}}}r_{{\mathcal {O}}_{i_{1}}}$, it is easy
to see that $w\lambda=d_{{\mathcal {O}}_{i_{1}}}d_{{\mathcal
{O}}_{i_{2}}}\cdots d_{{\mathcal {O}}_{i_{n}}}\lambda$ and
$\lambda$ satisfies $(\ast)$ if and only if $w\lambda$ satisfies
$(\ast)$.

\medskip
{\bf Corollary 3.4.} Suppose $G$ is an abelian group acting
on $Q$ admissibly. If $\mu=w\lambda$ for some $w\in
{\mathcal {W}}$ such that $w$ has minimal length with this property, then
there is an equivalence
$$\Pi_{Q}^{\lambda}G\mbox{-{\bf mod}}\longrightarrow
\Pi_{Q}^{\mu}G\mbox{-{\bf mod}},$$
which acts as $w$ on dimension vector.

\medskip
{\bf Proof.} By induction on the length of $w$, it deduces to
the case when $w$ is the reflection $r_{{\mathcal {O}}_{i}}$ for some
$i\in I$. Let $w=r_{{\mathcal {O}}_{i}}$, if $\lambda_{i}=0$, then
$w\lambda=\lambda$, this contradict to the minimality of $w$. So
$\lambda_{i}\neq 0$ and the corollary follows from Theorem 1.3.
\hfill$\square$

\section{Morita equivalence}
In this section, the aim is to prove Theorem 1.4.
Let us recall some definitions of group species and their
representations introduced in \cite{De}.

\medskip
{\bf Definition 4.1.} A group species is a triple $\Gamma=\big(I,
~(G_{i})_{i\in I}, ~(A_{i, j})_{(i, j)\in I^{2}}\big)$, where $I$ is
a set,  $G_{i}$ is a group for $i\in I$, and
$A_{i, j}$ is a $(kG_{j}, kG_{i})$-bimodule  for  $(i, j)\in
I^{2}$ ($kG_{j}$ acts on the left and $kG_{i}$ acts on the right).

\medskip
{\bf Definition 4.2.} A representation of $\Gamma$ is a pair
$\big((V_{i})_{i\in I}, ~(x_{i, j})_{(i, j)\in I^{2}}\big)$, where
 $V_{i}$ is a $k$-representation of $G_{i}$ for $i\in I$, and
 $x_{i, j}\in \Hom_{G_{i}}(A_{i,
j}\otimes_{G_{i}}V_{i}, ~V_{j})$ for $(i, j)\in I^{2}$.

\medskip
Let $V=\big((V_{i})_{i\in I}, ~(x_{i, j})_{(i, j)\in I^{2}}\big)$
and $V'=\big((V'_{i})_{i\in I}, ~(x'_{i, j})_{(i, j)\in I^{2}}\big)$
be two representations of $\Gamma$, the morphism from $V$ to $V'$ is
a family $ (f_{i})_{i\in I}\in \prod_{i\in I}\Hom_{G_{i}}(V_{i},
~V'_{i})$ such that the following diagram commute
\[\begin{array}{ccc}
A_{i, j}\otimes_{G_{i}}V_{i} & \stackrel{x_{i, j}}{\longrightarrow}
&
V_{j} \\
\llap{$\id_{A_{i, j}}\otimes f_{i}$}\downarrow  & & \downarrow
\rlap{$f_{j}$}\\
A_{i, j}\otimes_{G_{i}}V'_{i} & \stackrel{x'_{i,
j}}{\longrightarrow} & V'_{j}.
\end{array}\]
A representation $V=\big((V_{i})_{i\in I}, ~(x_{i, j})_{(i, j)\in
I^{2}}\big)$ is said to be finite if $\dim_{k}(V_{i})$ are all finite. It is easy
to check that all finite representations of a group species
$\Gamma$ together with the morphisms defined as above form an abelian category.
We denote it by $\Rep(\Gamma)$.

\medskip
Obviously, if $G_{i}$ is a trivial group for each $i\in I$, the category
$\Rep(\Gamma)$ coincide with the category of finite representations
of quiver.

Let $G$, $H$ be two groups. For any $(kG, kH)$-bimodule $V$ and
$(kH, kG)$-bimodule $V'$, a non degenerate paring between $V$ and
$V'$ is a bilinear map $\langle-,-\rangle$ from $V\times V'$ to $k$
such that for all $v\in V\setminus\{0\}$ and $v'\in
V'\setminus\{0\}$, $\langle v, -\rangle$ and $\langle-, v'\rangle$
do not vanish and such that for all $(g, h)\in G\times H$ and $(v,
v')\in V\times V'$, $\langle gvh,v'\rangle=\langle v,hv'g\rangle$.

\medskip
{\bf Definition 4.3.} The double group species $\Gamma^{d}$ is a
triple $\big(\Gamma, ~(\langle-,-\rangle_{i, j})_{(i, j)\in I^{2}},
~(\varepsilon_{i,j})_{(i, j)\in I^{2}}\big)$ where $\Gamma=\big(I,
~(G_{i})_{i\in I}, ~(A_{i, j})_{(i, j)\in I^{2}}\big)$ is a group
species and for each $(i, j)\in I^{2}$, $\langle-,-\rangle_{i, j}$
is a non degenerate paring between $A_{i, j}$ and $A_{j, i}$
satifying for each $(a, a')\in A_{i, j}\times A_{j, i}$, $\langle
a,a'\rangle_{i, j}=\langle a',a\rangle_{j, i}$, and
$\varepsilon_{i,j}$ is an automorphism of $A_{i, j}$ such that
$\varepsilon_{i,j}=-^{t}\varepsilon_{j,i}$ ($^{t}\varepsilon_{j,i}$
is the transpose of $\varepsilon_{j,i}$).

\medskip
By the definition, one see that $\Gamma^{d}$ is a double quiver
if $\Gamma$ is a quiver.

\medskip
Let now $\Gamma^{d}=\big(\Gamma, ~(\langle-,-\rangle_{i, j})_{(i,
j)\in I^{2}}, ~(\varepsilon_{i,j})_{(i, j)\in I^{2}}\big)$ be a
double group species. For each representation $\big((V_{i})_{i\in
I}, ~(x_{i, j})_{(i, j)\in I^{2}}\big)$ of $\Gamma^{d}$ and $(i,
j)\in I^{2}$ , there is an isomorphism \cite{De}:
$$
\begin{aligned}\varphi : \Hom_{G_{j}}(A_{i, j}
\otimes_{G_{i}}V_{i}, ~V_{j}) &\longrightarrow \Hom_{G_{i}}(V_{i},
~A_{j, i}
\otimes_{G_{j}}V_{j})\nonumber\\
f \qquad &\mapsto  \qquad\widetilde{f}
\end{aligned}
$$
where $\widetilde{f}(v)=\sum_{a\in {\mathcal {B}}_{i,
j}}a^{\ast}\otimes f(a\otimes v)$ for $v\in V_{i}$, ${\mathcal
{B}}_{i, j}$ is a basis of $A_{i, j}$ and $a^{\ast}$ is the element
of $A_{j, i}$ corresponding to $a$ in the dual basis (under the
bilinear form $\langle-,-\rangle_{i, j}$) of ${\mathcal {B}}_{i,
j}$. By this isomorphism, we give one definition as follows.

\medskip
{\bf Definition 4.4.} For any $\lambda=(\lambda_{i})_{i\in I} \in
k^{I}$, a representation $\big((V_{i})_{i\in I}, ~(x_{i,
j})_{(i,j)\in I^{2}}\big)$ of double group species
$\Gamma^{d}=\big(\Gamma, ~(\langle-,-\rangle_{i, j})_{(i, j)\in
I^{2}}, ~(\varepsilon_{i,j})_{(i, j)\in I^{2}}\big)$ is said to
satisfy the $\lambda$-preprojective relations if for each $i\in I$,
$$\sum_{j\in I}x_{j, i}\circ
(\varepsilon_{j, i}\otimes \id_{V_{j}})\circ \widetilde{x}_{i,
j}-\lambda_{i}=0.$$
Note that all representations of $\Gamma^{d}$
satisfying $\lambda$-preprojective relations form a full subcategory
of $\Rep(\Gamma)$, which is denoted by $\lambda\mbox{-}\Rep(\Gamma)$. In
particular, if $\Gamma$ is a double quiver, the
$\lambda$-preprojective relations coincide with the deformed
preprojective relations stated in Section 3.

\medskip
The following lemma is similar to \cite [Lemma 13]{De}. For the
completeness we skecth the proof.

\medskip
{\bf Lemma 4.1.} The category $\lambda\mbox{-}\Rep(\Gamma)$ does not
depend on the choice of $\varepsilon_{i,j}$ or
$\langle-,-\rangle_{i, j}$ up to isomorphism.

\medskip
{\bf Proof.}  It suffices to prove  there is a category
equivalence $$F : \Rep(\Gamma)\longrightarrow \Rep(\Gamma)$$ such
that $F$ sends the $\lambda$-preprojective relations of
$\Gamma^{d}_{1}=\big(\Gamma, ~(\langle-,-\rangle_{i, j})_{(i, j)\in
I^{2}}, ~(\varepsilon_{i,j})_{(i, j)\in I^{2}}\big)$ to the
$\lambda$-preprojective relations of $\Gamma^{d}_{2}=\big(\Gamma,
~(\langle-,-\rangle'_{i, j})_{(i, j)\in I^{2}},
~(\varepsilon'_{i,j})_{(i, j)\in I^{2}}\big)$, for any double group
species $\Gamma^{d}_{1}$ and $\Gamma^{d}_{2}$.

First of all, we prove that $\lambda\mbox{-}\Rep(\Gamma)$ does not
depend on the choice of $\varepsilon_{i,j}$. Given any $(i, j)\in
I^{2}$, we define an automorphism $\varphi_{i,j}$ of $A_{i,j}$ such
that $\varphi_{j,i}=\varepsilon'_{j,i}{^{t}\varphi_{i,j}}^{-1}
\varepsilon_{j,i}^{-1}$, and a functor
$$
\begin{aligned}F : \Rep(\Gamma)\qquad &\longrightarrow
\qquad\qquad\Rep(\Gamma)\nonumber\\
\big((V_{i})_{i\in I}, ~(x_{i, j})_{(i,j)\in I^{2}}\big)  &\mapsto
\big((V_{i})_{i\in I}, ~(x_{i, j}\circ(\varphi_{i,j}^{-1}\otimes
\id_{V_{i}}))_{(i,j)\in I^{2}}\big).
\end{aligned}
$$
One can check that $F$ is an equivalence of categories. $F$ sends
the representations of $\Gamma$ satisfying $\lambda$-preprojective
relations of type $\varepsilon_{i,j}$ to the representations of
$\Gamma$ satisfying $\lambda$-preprojective relations of type
$\varepsilon'_{i,j}$.

Secondly, suppose $\langle-,-\rangle'_{i, j}\neq
\langle-,-\rangle_{i, j}$. For any $(i, j)\in I^{2}$, any basis
${\mathcal {B}}_{i,j}$ of $A_{i,j}$, we denote by $d(x)$ and $d'(x)$
the dual basis of the element $x\in{\mathcal {B}}_{i,j}$ for the
pairing $\langle-,-\rangle_{i, j}$ and $\langle-,-\rangle'_{i, j}$
respectively. Let $\varepsilon''_{i,j}$ be the automorphism of
$A_{i,j}$ satisfying
$\varepsilon''_{i,j}(d'(x))=\varepsilon_{i,j}(d(x))$ for any
$x\in{\mathcal {B}}_{i,j}$. Note that $\varepsilon''_{i,j}$ does not
depend on the basis ${\mathcal {B}}_{i,j}$, we have
$\varepsilon''_{i,j}=-^{t}\varepsilon''_{j,i}\Leftrightarrow
\varepsilon''_{i,j}(d'(\varepsilon''_{j,i}(d'(x))))=-x, ~\forall
~x\in{\mathcal {B}}_{i,j}\Leftrightarrow
\varepsilon_{i,j}(d(\varepsilon_{j,i}(d(x))))=-x, ~\forall
~x\in{\mathcal {B}}_{i,j}\Leftrightarrow
\varepsilon_{i,j}=-^{t}\varepsilon_{j,i}$. Therefore, $\big(\Gamma,
~(\langle-,-\rangle'_{i, j})_{(i, j)\in I^{2}},
~(\varepsilon''_{i,j})_{(i, j)\in I^{2}}\big)$ is a double group
species and $$\sum_{j\in I}x_{j, i}\circ (\varepsilon_{j, i}\otimes
id_{V_{j}})\circ \widetilde{x}_{i, j}-\lambda_{i}=0 ~\mbox{ if and
only if } ~ \sum_{j\in I}x_{j, i}\circ (\varepsilon''_{j, i}\otimes
id_{V_{j}})\circ \widetilde{x}_{i, j}-\lambda_{i}=0,$$ for each
$i\in I$. Hence, It reduces the case that $\lambda\mbox{-}\Rep(\Gamma)$
does not depend on the choice of $\varepsilon_{i,j}$ again.
\hfill$\square$

\medskip
From now on, we view $\Lambda$ and $kG$ as
subalgebras of $\Lambda G$ and identify $\Lambda G$ with the algebra
$$\langle\Lambda, kG \mid g(a)=gag^{-1}, \forall ~g\in G,
\forall ~a\in \Lambda\rangle_{\mathrm{algebra}}$$ for any skew group
algebra $\Lambda G$.

Given a finite quiver $Q=(I, A)$ and a group $G\subseteq \Aut(Q)$, then
$G$ acts on the path algebra $kQ$ by permuting the set of primitive idempotents
$\{e_{i}\mid i\in I\}$.  It induces an action of $G$
on $k\overline{Q}$ in a natural way, i.e., $(a')^{g}=(a^{g})'$ for
any $a\in A$ and $g\in G$. As in Section 2, we get a skew group
algebra $\Pi_{Q}^{\lambda}G$, where $\lambda=(\lambda_{i})_{i\in I}$
satisfies $\lambda_{i}=\lambda_{j}$ if vertices $i$ and $j$ are in the
same $G$-orbit.

For each
$i\in I$, we denote by $G_{i}$ the subgroup of $G$ stabilizing $e_{i}$,
and by $$A_{i,j}:=e_{j}(\rad(kQ)/\rad^{2}(kQ))e_{i},$$ where
$\rad(kQ)$ is the Jacobson radical of $kQ$. Then we can get a group
species
$$\Gamma:=\big(I, ~(G_{i})_{i\in I}, ~(A_{i, j})_{(i, j)\in
I^{2}}\big)$$ and a double group species
$$\overline{\Gamma}^{d}=\big(\overline{\Gamma},
~(\langle-,-\rangle_{i, j})_{(i, j)\in I^{2}},
~(\varepsilon_{i,j})_{(i, j)\in I^{2}}\big),$$
where
$\overline{\Gamma}:=\big(I, ~(G_{i})_{i\in I}, ~(\overline{A}_{i,
j})_{(i, j)\in I^{2}}\big)$, $\overline{A}_{i, j}=A_{i, j}\oplus
A^{\ast}_{j, i}$, ~$A^{\ast}_{j, i}:= \Hom_{k}(A_{j, i}, ~k)$,
$\varepsilon_{i,j}=\id_{A_{i,j}}\oplus\left(-\id_{A^{\ast}_{j,i}}\right)$, and
$\langle-,-\rangle_{i, j}: \overline{A}_{i, j}\times
\overline{A}_{j, i}\rightarrow k$ given by
$$\langle a, ~b\rangle_{i,
j}=b_{2}(a_{1})+a_{2}(b_{1})$$ for $a=(a_{1},
a_{2})\in\overline{A}_{i, j}$ and $b=(b_{1}, b_{2})\in\overline{A}_{j,
i}$. It is easy to see that
$$
\begin{aligned}(1)& \quad\langle gag^{-1}, ~gbg^{-1}\rangle_{gi, gj}
=\langle a, b\rangle_{i, j}, \mbox {~~for~~} a\in \overline{A}_{i,
j}, ~b\in \overline{A}_{j,
i} \mbox{~~and~~} g\in G; \qquad\qquad\qquad\nonumber\\
(2)& \quad\varepsilon_{gi,gj}(gag^{-1})=g\varepsilon_{i,j}(a)g^{-1},
\mbox {~~for~~} a\in \overline{A}_{i, j}, \mbox{~~and~~} g\in G
\end{aligned}
$$ where we denote $gi:=g(i)$ for convenience.

Let $\widetilde{I}$ be a set of representatives of the classes of
$I$ under the action of $G$ on $I$, $E_{i, j}$ a set of
representatives of the classes of $\widetilde{I}^{2}$ under the
action of $G$ on ${\mathcal {O}}_{i}\times {\mathcal {O}}_{j}$,
where ${\mathcal {O}}_{i}$ is the $G$-orbit of $i$. If $i_{0}\in
\widetilde{I}$ is a representative of the class of $i$, we denote
by $\kappa_{i}$ the element of $G$ such that $\kappa_{i}i_{0}=i$.

\medskip
Now, we want define a new double group species. Firstly, for any
$(i, j)\in \widetilde{I}^{2}$, we let
$$\widetilde{A}_{i, j}=\bigoplus_{(i', j')\in E_{i, j}}G_{j}
\kappa^{-1}_{j'}\overline{A}_{i', j'}\kappa_{i'}G_{i}.$$
It is easy to check that $\widetilde{A}_{i, j}$ is a $(kG_{j},
kG_{i})$-bimodule and $\widetilde{A}_{i, j}$ does not depend on the
choice of $E_{i, j}$. Therefore, we get a group species
$\widetilde{\Gamma}_{G}:=\big(\widetilde{I}, ~(G_{i})_{i\in
\widetilde{I}}, ~(\widetilde{A}_{i, j})_{(i, j)\in
\widetilde{I}^{2}}\big)$. Secondly, for any $(i, j)\in
\widetilde{I}^{2}$, $h_{1} \kappa^{-1}_{j'}a\kappa_{i'}g_{1}\in
\widetilde{A}_{i, j}$, $g_{2} \kappa^{-1}_{i''}b\kappa_{j''}h_{2}\in
\widetilde{A}_{j, i}$, we define
$$\begin{aligned}\langle h_{1}
\kappa^{-1}_{j'}a\kappa_{i'}g_{1},& ~g_{2}
\kappa^{-1}_{i''}b\kappa_{j''}h_{2}\rangle'_{i, j}\nonumber\\
=& {\left\{\begin{array}{ll} \langle \kappa_{j''}h_{2}h_{1}
\kappa^{-1}_{j'}a\kappa_{i'}g_{1}g_{2}\kappa^{-1}_{i''},
~b\rangle_{i', j'}, & \mbox{if }(i', j')=(i'', j''), \\&\mbox{ and
}\kappa_{j''}h_{2}h_{1}
\kappa^{-1}_{j'}\kappa_{i'}g_{1}g_{2}\kappa^{-1}_{i''}=e;\\
0,& \mbox{otherwise}.
\end{array}\right.}
\end{aligned}$$
and $\varepsilon'_{i,j}(h_{1}
\kappa^{-1}_{j'}a\kappa_{i'}g_{1})=\frac{1}{|G_{i}||G_{j}|}
h_{1}\kappa^{-1}_{j'}\varepsilon_{i,j}(a)\kappa_{i'}g_{1}$. Then,
$\big(\widetilde{\Gamma}_{G}, ~(\langle-,-\rangle'_{i, j})_{(i,
j)\in I^{2}}, ~(\varepsilon'_{i,j})_{(i, j)\in I^{2}}\big)$ is a
double group species (see \cite [Lemma 16]{De}). We define
the functor
$$
\begin{aligned}\Phi :  \quad(k\overline{Q})G\mbox{-}\Mod\qquad &\longrightarrow
 \qquad\qquad\Rep(\widetilde{\Gamma}_{G}) \nonumber\\
 V \qquad\qquad &\mapsto
\qquad \big((W_{i})_{i\in \widetilde{I}}, ~(y_{i, j})_{(i, j)\in
\widetilde{I}^{2}}\big)
\end{aligned}
$$
where $W_{i}=e_{i}V$ for each $i\in \widetilde{I}$ and $y_{i, j}:
\widetilde{A}_{i, j}\otimes_{G_{i}}W_{i}\rightarrow W_{j}$, $y_{i,
j}(a\otimes v)=av$ for each $(i, j)\in \widetilde{I}^{2}$. Then, by
\cite [Proposition 15]{De} we have

{\bf Lemma 4.2.} $\Phi:(k\overline{Q})G\mbox{-}\Mod \rightarrow
\Rep(\widetilde{\Gamma}_{G})$ is a categories equivalence.
\hfill$\square$

For $i, j\in I$, we define an equivalence relation
$``\sim"$ on $G_{i}\times G_{j}$ by $(g, h) \sim (g', h')$ if and
only if $h\kappa_{j}^{-1}\kappa_{i}g=h'\kappa_{j}^{-1}\kappa_{i}g'$,
and denote by $G_{i, j}$ the set of representatives of the
equivalence classes for this relation, then
$$\widetilde{{\mathcal {B}}}_{i,j}=\bigcup_{(i', j')\in E_{i, j}}
\bigcup_{(g, h)\in G_{i', j'}}h\kappa_{j'}^{-1}\overline{{\mathcal
{B}}}_{i', j'}\kappa_{i'}g$$ is a basis of $\widetilde{A}_{i, j}$
for any $(i, j)\in I^{2}$, where $\overline{{\mathcal {B}}}_{i,j}$
is a basis of $\overline{A}_{i, j}$ \cite{De}. One can check that,
if $h_{1} \kappa^{-1}_{j'}a\kappa_{i'}g_{1}\in \widetilde{{\mathcal
{B}}}_{i,j}$, the dual basis of
$\widetilde{A}_{i,j}$ in $\widetilde{{\mathcal
{B}}}_{i,j}$ for $\langle-,-\rangle'_{ij}$ is $(h_{1}
\kappa^{-1}_{j'}a\kappa_{i'}g_{1})^{\ast}=g_{1}^{-1}\kappa_{i'}^{-1}
a^{\ast}\kappa_{j'}h^{-1}_{1}$, where $a^{\ast}$ is the dual basis of
$A_{j',i'}$ corresponding to $a$ in ${\mathcal
{B}}_{i',j'}$ for $\langle-,-\rangle_{i',j'}$.

As the proof of \cite [Proposition. 17]{De}, we have

{\bf Proposition 4.3.} If $\lambda=(\lambda_{i})_{i\in I}$ satisfies
$\lambda_{i}=\lambda_{j}$ for any vertices $i$, $j$ in the same
$G$-orbit, then the category of representations of
$(k\overline{Q})G$ satisfying the $\lambda$-preprojective relations
is equivalent to the category of representations of
$\widetilde{\Gamma}_{G}$ satisfying the $\xi$-preprojective
relations, where $\xi=(\xi_{i_{0}})_{i_{0}\in \widetilde{I}}$,
$\xi_{i_{0}}=\frac{\lambda_{i}}{|G_{i}|}$ if $i\in {\mathcal
{O}}_{i_{0}}$.

{\bf Proof.} Note that $|G_{i}|=|G_{j}|$ if $i$ and $j$ are
in the same $G$-orbit. The hypothesis of
$\lambda=(\lambda_{i})_{i\in I}$, $\xi=(\xi_{i_{0}})_{i_{0}\in
\widetilde{I}}$ are given well.

Recall that there is an equivalence functor
$\Phi:(k\overline{Q})G\mbox{-}\Mod \rightarrow
\Rep(\widetilde{\Gamma}_{G})$ and every $(k\overline{Q}) G$
representation $V=\big((V_{i})_{i\in I}, ~(x_{i, j})_{(i, j)\in
I^{2}}\big)$ is $G$-invariant. It  follows that
$$\sum_{j\in I}x_{j, i}\circ (\varepsilon_{j, i}\otimes \id_{V_{j}})
\circ \widetilde{x}_{i,j}-\lambda_{i}\id_{V_{i}}=0$$
for any $i\in I$ is equivalent to
$\sum_{j\in I}x_{j, i}\circ (\varepsilon_{j, i}\otimes
\id_{V_{j}})\circ \widetilde{x}_{i, j}-\lambda_{i}\id_{V_{i}}=0$,
for any $i\in \widetilde{I}$. For this reason, it reduces to prove that
$\Phi(V)=\big((W_{i})_{i\in \widetilde{I}}, ~(y_{i, j})_{(i, j)\in
\widetilde{I}^{2}}\big)$ satisfies $\sum_{j\in I}y_{j, i}\circ
(\varepsilon'_{j, i}\otimes \id_{W_{j}})\circ \widetilde{y}_{i,
j}-\xi_{i}\id_{W_{i}}=0$ if and only if $V=\big((V_{i})_{i\in I},
~(x_{i, j})_{(i, j)\in I^{2}}\big)$ satisfies $\sum_{j\in I}x_{j,
i}\circ (\varepsilon_{j, i}\otimes \id_{V_{j}})\circ
\widetilde{x}_{i, j}-\lambda_{i}\id_{V_{i}}=0$ for any $i\in
\widetilde{I}$ and $(k\overline{Q}) G$-module $V$.
For all $i\in \widetilde{I}$ and $v\in W_{i}$, we have
$$\begin{aligned}
\lambda_{i}v=& \sum_{j\in \widetilde{I}}y_{j, i}\circ
(\varepsilon'_{j, i}\otimes id_{W_{j}})\circ \widetilde{y}_{i, j}(v)
=\sum_{j\in \widetilde{I}}\sum_{a\in \widetilde{{\mathcal
{B}}}_{i, j}}\varepsilon'_{j, i}(a^{\ast}) av\nonumber\\
=&\sum_{j\in \widetilde{I}}\sum_{(i', j')\in E_{i, j}}\sum_{(g,
h)\in G_{i', j'}}\sum_{a\in {\mathcal {B}}_{i', j'}}
\varepsilon'_{j,
i}((g^{-1}\kappa_{i'}^{-1}a\kappa_{j'}h^{-1})^{\ast})
h\kappa_{j'}^{-1}a\kappa_{i'}gv\nonumber\\
=&\sum_{j\in \widetilde{I}}\sum_{(i', j')\in E_{i, j}}\sum_{(g,
h)\in G_{i', j'}}\sum_{a\in {\mathcal {B}}_{i',
j'}}\frac{1}{|G_{i'}||G_{j'}|}g^{-1}\kappa_{i'}^{-1}\varepsilon_{j',
i'}(a^{\ast})
a\kappa_{i'}gv\nonumber\\
=&\frac{1}{|G_{i}|}\sum_{j\in
\widetilde{I}}\frac{1}{|G_{j}|}\sum_{(i', j')\in E_{i, j}}\sum_{(g,
h)\in G_{i', j'}}\sum_{a\in {\mathcal {B}}_{i',
j'}}\varepsilon_{g^{-1}\kappa_{i'}^{-1}j',i}(g^{-1}\kappa_{i'}^{-1}a^{\ast}\kappa_{i'}g)g^{-1}\kappa_{i'}^{-1}
a\kappa_{i'}gv\nonumber\\
=&\frac{1}{|G_{i}|}\sum_{j\in
\widetilde{I}}\frac{1}{|G_{j}|}\sum_{(i', j')\in E_{i, j}}\sum_{(g,
h)\in G_{i', j'}}\sum_{a\in {\mathcal {B}}_{i,
g^{-1}\kappa_{i'}^{-1}j'}}\varepsilon_{g^{-1}
\kappa_{i'}^{-1}j',i}(a^{\ast})av\nonumber\\
=&\frac{1}{|G_{i}|}\sum_{j\in
\widetilde{I}}\frac{1}{|G_{j}|}\sum_{(i', j')\in E_{i, j}}\sum_{(g,
h)\in G_{i', j'}}x_{g^{-1}\kappa_{i'}^{-1}j',
i}\circ(\varepsilon_{g^{-1}\kappa^{-1}_{i'}j', i}\otimes
\id_{W_{g^{-1}\kappa^{-1}_{i'}j'}})\circ\widetilde{x}_{i,
g^{-1}\kappa_{i'}^{-1}j'}(v).
\end{aligned}$$

Define a map $\zeta_{i', j'} : G_{i', j'}\rightarrow I$ by
$(g, h)\mapsto g^{-1}\kappa^{-1}_{i'}j'$ for all $(i', j')\in E_{i,
j}$, $i, j \in I$. It is easy to check that $\Im \zeta_{i', j'}=I_{i',
j'}:=\{j''\in I\mid \exists ~g\in G, \mbox{ such that } (i,
j'')=g(i', j')\}$ and $|\zeta_{i', j'}^{-1}(j'')|=|G_{j}|$ for any
$j''\in I_{i', j'}$, where $\zeta_{i', j'}^{-1}(j'')$ is the set of
preimages of $j''$. Note that
$$I=\bigcup_{j\in \widetilde{I}}\bigcup_{(i', j')\in E_{i,
j}}I_{i', j'},$$
it follows that
$$\begin{aligned}
&\sum_{j\in \widetilde{I}}\frac{1}{|G_{j}|}\sum_{(i', j')\in E_{i,
j}}\sum_{(g, h)\in G_{i', j'}}x_{g^{-1}\kappa_{i'}^{-1}j',
i}\circ(\varepsilon_{g^{-1}\kappa^{-1}_{i'}j', i}\otimes
\id_{W_{g^{-1}\kappa^{-1}_{i'}j'}})\circ\widetilde{x}_{i,
g^{-1}\kappa_{i'}^{-1}j'}(v)-\xi_{i}v\nonumber\\
&=\sum_{j\in I}x_{j, i}\circ(\varepsilon_{j, i}\otimes
id_{V_{j}})\circ\widetilde{x}_{i, j}(v)-\lambda_{i}v
\end{aligned}$$
for all $i\in \widetilde{I}$. The proof is completed.
\hfill$\square$

In the following we define a new quiver $Q_{G}=\big(I_{G}, ~(A'_{(i, \rho), (j,
\sigma)})_{((i, \rho), (j, \sigma))\in I^{2}_{G}}\big)$ associated
with the group species $\widetilde{\Gamma}_{G}=\big(\widetilde{I},
~(G_{i})_{i\in \widetilde{I}}, ~(\widetilde{A}_{i, j})_{(i, j)\in
\widetilde{I}^{2}}\big)$ by
$$
I_{G}=\bigcup_{i\in \widetilde{I}}\{i\}\times \irr(G_{i}) \quad
\mbox{ and } \quad A'_{(i, \rho), (j,
\sigma)}=\Hom_{G_{j}}(\widetilde{A}_{i, j}\otimes_{G_{i}}\rho,
~\sigma)^{\ast},
$$
where $(i, \rho), (j, \sigma)\in I_{G}$,
$\Hom_{G_{j}}(\widetilde{A}_{i, j}\otimes_{G_{i}}\rho,
~\sigma)^{\ast}=\Hom_{k}\big(\Hom_{G_{j}}(\widetilde{A}_{i,
j}\otimes_{G_{i}}\rho, ~\sigma), ~k\big)$, and $\irr(G_{i})$ is the
set of representatives of isomorphism classes of irreducible
representations of $G_{i}$.

For $(i, \rho), (j, \sigma)\in I_{G}$ and two vector species $V$, $W$,
there is a linear isomorphism $$\varphi_{(i, \rho), (j, \sigma)}:
\Hom_{k}\big(\Hom_{G_{j}}(\widetilde{A}_{i,j}\otimes_{G_{i}}\rho,
\sigma)^{\ast}\otimes_{k}V, ~W\big)\rightarrow
\Hom_{G_{j}}(\widetilde{A}_{i,j}\otimes_{G_{i}}\rho\otimes_{k}V,
~\sigma\otimes_{k}W)$$ define by $\varphi_{(i, \rho), (j,
\sigma)}(f)(a\otimes t\otimes v)=\sum_{g\in {\mathcal {B}}_{(i,
\rho), (j, \sigma)}}g(a\otimes t)\otimes f(g^{\ast}\otimes v)$,
which does not depend on the choice of a basis ${\mathcal {B}}_{(i,
\rho), (j, \sigma)}$ of
$\Hom_{G_{j}}(\widetilde{A}_{i,j}\otimes_{G_{i}}\rho, \sigma)$.

By \cite [Proposition 19]{De}, if we define the functor
$$
\begin{aligned}\Psi : \quad kQ_{G}\mbox{-}\Mod
\qquad\qquad\qquad &\longrightarrow
 \qquad\qquad \Rep(\widetilde{\Gamma}_{G})\nonumber\\
\big((V_{(i,\rho)})_{(i,\rho)\in I_{G}}, ~(x_{(i, \rho), (j,
\sigma)})_{((i, \rho), (j, \sigma))\in I_{G}^{2}}\big) &\mapsto
\qquad \big((W_{i})_{i\in \widetilde{I}}, ~(y_{i, j})_{(i, j)\in
\widetilde{I}^{2}}\big)
\end{aligned}
$$
where $W_{i}=\bigoplus_{\rho\in
\irr(G_{i})}\rho\otimes_{k}V_{(i,\rho)}$ and $y_{i,
j}=\Big(\varphi_{(i, \rho), (j, \sigma)}(x_{(i, \rho), (j,
\sigma)})\Big)_{|\irr(G_{j})|\times|\irr(G_{i})|}$ for $i, j\in
\widetilde{I}$ .
We have

\medskip
{\bf Lemma 4.4.} $\Psi : kQ_{G}\mbox{-}\Mod \rightarrow
\Rep(\widetilde{\Gamma}_{G})$ is a categories equivalence.
\hfill$\square$

We claim that $Q_{G}$ is a double quiver. Indeed, let $(i, \rho), (j,
\sigma)\in I_{G}$, $f\in A'^{\ast}_{(i, \rho), (j,
\sigma)}=\Hom_{G_{j}}(\widetilde{A}_{i, j}\otimes_{G_{i}}\rho,
~\sigma)$, and $g\in A'^{\ast}_{(j, \sigma), (i,
\rho)}=\Hom_{G_{i}}(\widetilde{A}_{j, i}\otimes_{G_{j}}\sigma,
~\rho)$, we define a map
$$
\begin{aligned}\phi_{a}: \rho \qquad &\longrightarrow
\qquad\rho \nonumber\\
 r \qquad &\mapsto g(a^{\ast}\otimes f(a\otimes r))
\end{aligned}
$$
where $a\in \widetilde{{\mathcal {B}}}_{i, j}$, $a^{\ast}$ is the
element of $\widetilde{A}_{j, i}$ corresponding to $a$ in the dual
basis of $\widetilde{{\mathcal {B}}}_{i, j}$ for
$\langle-,-\rangle'_{i,j}$. It is easy to see that $\phi_{a}$ is a
$G_{i}$-representation homomorphism. Since $\rho$ is irreducible, we
know that $\phi_{a}$ is a multiple of the identity and
$g(a^{\ast}\otimes f(a\otimes
r))=\frac{1}{\dim_{k}\rho}\Tr\phi_{a}\cdot r$, where $\Tr\phi_{a}$
is the trace of $\phi_{a}$. Now, we let
$$
\langle f, g\rangle_{(i, \rho), (j, \sigma)}=\sum_{a\in
\widetilde{{\mathcal {B}}}_{i,j}}\Tr\phi_{a}\quad \mbox{ and } \quad
\varepsilon_{(j, \sigma), (i,
\rho)}(f)=f\circ({^{t}\varepsilon'_{j,i}}\otimes
\id_{\rho})=-f\circ(\varepsilon'_{i,j}\otimes \id_{\rho}),
$$
for $(i, \rho), (j, \sigma)\in I_{G}$. One can check that $\langle
-, -\rangle_{(i, \rho), (j, \sigma)}$ is a non degenerate paring and
$\varepsilon_{(j, \sigma), (i, \rho)}=-{^{t}\varepsilon_{(i, \rho),
(j, \sigma)}}$. Thus, one has also a non degenerate
paring between $A'_{(i, \rho), (j, \sigma)}$ and $A'_{(j, \sigma),
(i, \rho)}$, which we also denote by $\langle -, -\rangle_{(i,
\rho), (j, \sigma)}$. The automorphism $\varepsilon_{(i, \rho), (j,
\sigma)}$ can be seen as an automorphism of $A'_{(i, \rho), (j,
\sigma)}$ through $\langle -, -\rangle_{(i, \rho), (j, \sigma)}$,
and therefore $Q_{G}$ is a double quiver.

\medskip
{\bf Proposition 4.5.} For any $\xi=(\xi_{i})_{i\in \widetilde{I}}$,
the category of representations of $\widetilde{\Gamma}_{G}$
satisfying the $\xi$-preprojective relations is equivalent to the
category of representations of $Q_{G}$ satisfying the
$\eta$-preprojective relations, where $\eta=(\eta_{(i, \rho)})_{(i,
\rho)\in I_{G}}$ satisfying $\eta_{(i,
\rho)}=\dim_{k}\rho\cdot\xi_{i}$.

\medskip
{\bf Proof.} Denote by ${\mathcal {B}}_{(i, \rho), (j, \sigma)}$ a
basis of $\Hom_{G_{j}}(\widetilde{A}_{i,j}\otimes_{G_{i}}\rho,
\sigma)$ such that ${\mathcal {B}}_{(i, \rho), (j, \sigma)}$ and
${\mathcal {B}}_{(j, \sigma), (i, \rho)}$ form a dual basis for
$\langle-,-\rangle_{(j, \sigma), (i, \rho)}$. Let $\alpha\in
{\mathcal {B}}_{(i, \rho), (j, \sigma)}$, $\beta\in{\mathcal
{B}}_{(j, \sigma), (i, \tau)}$, we define map
$\varphi_{\alpha,\beta}:\rho\rightarrow \tau$, by
$\varphi_{\alpha,\beta}(r)=\sum_{a\in \widetilde{{\mathcal
{B}}}_{i,j}}\beta(\varepsilon'_{j,i}(a^{\ast})\otimes\alpha(a\otimes
r))$ for any $r\in\rho$. One can check that $\varphi_{\alpha,\beta}$
does not depend on the choice of a basis $\widetilde{{\mathcal
{B}}}_{i,j}$ of $\widetilde{A}_{i,j}$, and $\varphi_{\alpha,\beta}$
is a $G_{i}$-representations homomorphism (or see \cite [Proposition
20]{De}).

Therefore, $\varphi_{\alpha,\beta}=0$ if $\tau\neq\rho$, since
$\tau$ and $\rho$ are irreducible. And if $\tau=\rho$, for any
$r\in\rho$, there is
$$
\begin{aligned}\varphi_{\alpha,\varepsilon^{-1}_{(i, \rho),
(j, \sigma)}(\beta)}(r)=&\sum_{a\in\widetilde{{\mathcal
{B}}}_{i,j}}\varepsilon^{-1}_{(i, \rho), (j,
\sigma)}(\beta)(\varepsilon'_{j,i}(a^{\ast})\otimes\alpha(a\otimes
r))\nonumber\\
=&-\sum_{a\in\widetilde{{\mathcal
{B}}}_{i,j}}(\beta\circ(\varepsilon'^{-1}_{j,i}\otimes\id_{\sigma}))
(\varepsilon'_{j,i}(a^{\ast})\otimes\alpha(a\otimes r))\nonumber\\
=&-\sum_{a\in\widetilde{{\mathcal
{B}}}_{i,j}}\beta(a^{\ast}\otimes\alpha(a\otimes
r))=-\frac{\langle\alpha,\beta\rangle_{(i, \rho), (j,
\sigma)}}{\dim_{k}\rho}r
\end{aligned}
$$

Now, given any $Q_{G}$-representation
$V=\big((V_{(i,\rho)})_{(i,\rho)\in I_{G}}, ~(x_{(i, \rho), (j,
\sigma)})_{((i, \rho), (j, \sigma))\in I_{G}^{2}}\big)$, if
$\Psi(V)=\big((W_{i})_{i\in \widetilde{I}}, ~(y_{i, j})_{(i, j)\in
\widetilde{I}^{2}}\big)$ satisfies the $\xi$-preprojective
relations, that is,  for all $i\in \widetilde{I}$, $r\otimes v\in
\rho\otimes V_{(i,\rho)}\subseteq W_{i}$,
$$
\begin{aligned}&\xi_{i}(r\otimes v)\nonumber\\
=&\sum_{i\in\widetilde{I}}y_{j,i}\circ
(\varepsilon'_{j,i}\otimes\id_{W_{j}})\circ
\widetilde{y}_{i,j}(r\otimes v)=\sum_{j\in\widetilde{I}}y_{j,i}
\big((\varepsilon'_{j,i}\otimes\id_{W_{j}})
\big(\sum_{a\in\widetilde{{\mathcal {B}}}_{i,j}}a^{\ast}\otimes
y_{i,j}(a\otimes r\otimes
v)\big)\big)\nonumber\\
=&\sum_{j\in\widetilde{I}}y_{j,i}\big((\varepsilon'_{j,i}
\otimes\id_{W_{j}})\big(\sum_{a\in\widetilde{{\mathcal
{B}}}_{i,j}}a^{\ast}\otimes\sum_{\sigma\in \irr(G_{j})}
\sum_{a\in{\mathcal {B}}_{(i, \rho), (j, \sigma)}} \alpha(a\otimes
r)\otimes x_{(i, \rho), (j, \sigma)}
(\alpha^{\ast}\otimes v)\big)\big)\nonumber\\
=&\sum_{j\in\widetilde{I}}\sum_{a\in\widetilde{{\mathcal
{B}}}_{i,j}}\sum_{\sigma\in \irr(G_{j})}\sum_{a\in{\mathcal
{B}}_{(i, \rho), (j,
\sigma)}}y_{j,i}\big(\varepsilon'_{j,i}(a^{\ast})\otimes\alpha(a\otimes
r)\otimes x_{(i, \rho), (j, \sigma)}(\alpha^{\ast}\otimes
v)\big)\nonumber\\
=&\sum_{j\in\widetilde{I}}\sum_{\alpha\in\widetilde{{\mathcal
{B}}}_{i,j}}\sum_{\substack{\sigma\in \irr(G_{j})\\
\tau\in \irr(G_{i})}}\sum_{\substack{\alpha\in{\mathcal {B}}_{(i,
\rho), (j, \sigma)}\\\beta\in{\mathcal {B}}_{(j, \sigma), (i,
\tau)}}}\beta\big(\varepsilon'_{j,i}(a^{\ast})\otimes\alpha(a\otimes
r)\big)\otimes x_{(j, \sigma), (i, \rho)}\big(\beta^{\ast}\otimes
x_{(i, \rho), (j, \sigma)}(\alpha^{\ast}\otimes v)\big)\nonumber\\
=&\sum_{(j,\sigma)\in I_{G}}\sum_{\substack{\alpha\in{\mathcal
{B}}_{(i, \rho), (j, \sigma)}\\\beta\in{\mathcal {B}}_{(j, \sigma),
(i, \tau)}}}\varphi_{\alpha,\beta}(r)\otimes x_{(j, \sigma), (i,
\rho)}\big(\beta^{\ast}\otimes x_{(i, \rho), (j,
\sigma)}(\alpha^{\ast}\otimes v)\big)\nonumber\\
=&\sum_{(j,\sigma)\in I_{G}}\sum_{\substack{\alpha\in{\mathcal
{B}}_{(i, \rho), (j, \sigma)}\\\beta\in{\mathcal {B}}_{(j, \sigma),
(i, \tau)}}}\varphi_{\alpha,\varepsilon^{-1}_{(i, \rho), (j,
\sigma)}(\beta)}(r)\otimes x_{(j, \sigma), (i,
\rho)}\big((\varepsilon^{-1}_{(i, \rho), (j,
\sigma)}(\beta))^{\ast}\otimes x_{(i, \rho), (j,
\sigma)}(\alpha^{\ast}\otimes v)\big)\nonumber\\
=&-\sum_{(j,\sigma)\in I_{G}}\sum_{\substack{\alpha\in{\mathcal
{B}}_{(i, \rho), (j, \sigma)}\\\beta\in{\mathcal {B}}_{(j, \sigma),
(i, \tau)}}}\frac{\langle\alpha,\beta\rangle_{(i, \rho), (j,
\sigma)}}{\dim_{k}\rho}r\otimes x_{(j, \sigma), (i,
\rho)}\big({^{t}\varepsilon}_{(i, \rho), (j,
\sigma)}(\beta^{\ast})\otimes x_{(i, \rho), (j,
\sigma)}(\alpha^{\ast}\otimes v)\big)\nonumber\\
=&\frac{r}{\dim_{k}\rho}\otimes\sum_{(j,\sigma)\in
I_{G}}\sum_{\alpha\in{\mathcal {B}}_{(i, \rho), (j, \sigma)}}x_{(j,
\sigma), (i, \rho)}\big(\varepsilon_{(j, \sigma), (i,
\rho)}(\alpha)\otimes x_{(i, \rho), (j,
\sigma)}(\alpha^{\ast}\otimes v)\big)\nonumber\\
=&\frac{r}{\dim_{k}\rho}\otimes\sum_{(j,\sigma)\in I_{G}}x_{(j,
\sigma), (i, \rho)}\circ(\varepsilon_{(j, \sigma), (i,
\rho)}\otimes\id_{V_{(j, \sigma)}})\circ \widetilde{x}_{(i, \rho),
(j, \sigma)}(v).
\end{aligned}
$$
This is equivalent to say that
$$\sum_{(j,\sigma)\in I_{G}}x_{(j,
\sigma), (i, \rho)}\circ(\varepsilon_{(j, \sigma), (i,
\rho)}\otimes\id_{V_{(j, \sigma)}})\circ \widetilde{x}_{(i, \rho),
(j, \sigma)}(v)-\eta_{(i,\rho)}(v)=0,$$ for any $v\in V_{(i,\rho)}$.
Hence we get the proof. \hfill$\square$

\medskip
{\bf Proof of Theorem 1.4.}
It follows from Proposition 4.3 and Proposition 4.5.
\hfill$\square$

By the definition,  for any quiver $Q$ we can get the quiver $Q_{G}$,
given by the double group species corresponding to $Q$.

In particular, if $Q$ is a Dynkin quiver, we can describe the quiver $Q_{G}$
immediately by direct calculation, which are given as the following table.

$$\begin{tabular}{|c|c|c|} \hline $Q$ & $G$ & $Q_{G}$ \\ \hline
$\setlength{\unitlength}{1mm} \begin {picture}(45,20)
\put(0,16){$A_{2n+1} (n\geq 0):$}\put(2,11){$1$}\put(2,1){$1'$}
\put(4,2){\vector (1,0){5}}\put(4,12){\vector (1,0){5}}
\put(10,11){$2$}\put(10,1){$2'$} \put(12,2){\vector
(1,0){5}}\put(12,12){\vector (1,0){5}}
\put(18,1.5){$\cdots$}\put(18,11.5){$\cdots$} \put(23,2){\vector
(1,0){5}}\put(23,12){\vector (1,0){5}}
\put(29,11){$n$}\put(29,1){$n'$} \put(32,3){\vector
(3,2){5}}\put(32,11){\vector (3,-2){5}} \put(38,6){$n+1$}\end
{picture}$ & $\setlength{\unitlength}{1mm}\begin {picture}(8,20)
\put(3,8){${\mathbb{Z}}_{2}$}\end {picture}$ &
$\setlength{\unitlength}{1mm} \begin {picture}(50,16)
\put(27,5){$n$}\put(30,7){\vector (3,2){5}}
\put(37,10){$n+1'$}\put(30,6){\vector (3,-2){5}}
\put(37,1){$n+1$}\put(0,5){$1$} \put(2,6){\vector
(1,0){5}}\put(8,5){$2$} \put(10,6){\vector
(1,0){5}}\put(16,5.5){$\cdots$} \put(21,6){\vector (1,0){5}}\end
{picture}$ \\
\hline $\setlength{\unitlength}{1mm} \begin {picture}(45,16)
\put(0,12){$D_{4}:$}\put(22,5){$1$}\put(24,7){\vector (3,2){7}}
\put(33,11){$2''$}\put(24,6){\vector (3,-2){7}}
\put(33,0){$2'$}\put(21,6){\vector (-1,0){7}} \put(11.5,5){$2$}\end
{picture}$ & $\setlength{\unitlength}{1mm}\begin {picture}(8,20)
\put(3,8){$S_{3}$}\end {picture}$ & $\setlength{\unitlength}{1mm}
\begin {picture}(45,16)
\put(10,14){$1'$}\put(10,7){1}\put(10,0){$1''$}\put(13,15){\vector
(1,0){17}}\put(13,1){\vector (1,0){17}}\put(31,14){2}\put(31,0){$2'$}
\put(13,8){\vector (3,-1){16}}\put(13,8.5){\vector
(3,1){16}}\end{picture}$ \\
\hline $\setlength{\unitlength}{1mm} \begin {picture}(50,16)
\put(0,12){$D_{n} (n\geq 0):$}\put(27,5){$n-2$}\put(36,7){\vector
(3,2){5}} \put(42,10){$n-1'$}\put(36,6){\vector (3,-2){5}}
\put(42,1){$n-1$}\put(0,5){$1$} \put(2,6){\vector
(1,0){5}}\put(8,5){$2$} \put(10,6){\vector
(1,0){5}}\put(16,5.5){$\cdots$} \put(21,6){\vector (1,0){5}}\end
{picture}$ & $\setlength{\unitlength}{1mm}\begin {picture}(8,20)
\put(3,8){${\mathbb{Z}}_{2}$}\end {picture}$ &
$\setlength{\unitlength}{1mm} \begin {picture}(52,20)
\put(2,11){$1$}\put(2,1){$1'$} \put(4,2){\vector
(1,0){5}}\put(4,12){\vector (1,0){5}}
\put(10,11){$2$}\put(10,1){$2'$} \put(12,2){\vector
(1,0){5}}\put(12,12){\vector (1,0){5}}
\put(18,1.5){$\cdots$}\put(18,11.5){$\cdots$} \put(23,2){\vector
(1,0){5}}\put(23,12){\vector (1,0){5}}
\put(29,11){$n-2$}\put(29,1){$n-2'$} \put(38,3){\vector
(3,2){5}}\put(38,11){\vector (3,-2){5}} \put(44,6){$n-1$}\end
{picture}$ \\
\hline $\setlength{\unitlength}{1mm} \begin {picture}(45,16)
\put(0,12){$E_{6}:$}\put(22,11){$1$}\put(22.5,10){\vector
(0,-1){5}}\put(22,1){$2$}\put(25,2){\vector (1,0){5}}
\put(31,1){$3$}\put(34,2){\vector (1,0){5}}\put(40,1){$4$}
\put(21,2){\vector (-1,0){5}}\put(13,1){$3'$} \put(12,2){\vector
(-1,0){5}}\put(4,1){$4'$}
\end {picture}$ & $\setlength{\unitlength}{1mm}\begin {picture}(8,20)
\put(3,8){${\mathbb{Z}}_{2}$}\end {picture}$ &
$\setlength{\unitlength}{1mm} \begin {picture}(45,16)
\put(22,11){$4$}\put(22.5,5){\vector
(0,1){5}}\put(22,1){$3$}\put(30,2){\vector (-1,0){5}}
\put(31,1){$2$}\put(39,2){\vector (-1,0){5}}\put(40,1){$1$}
\put(16,2){\vector (1,0){5}}\put(13,1){$2'$} \put(7,2){\vector
(1,0){5}}\put(4,1){$1'$}
\end {picture}$ \\
\hline
\end{tabular}$$
Since the deformed preprojective algebra $\Pi_{Q}^{\lambda}$
does not depend on the orientation of $Q$, the above table shows that for
all the
Dynkin quivers with relevant nontrivial automorphism group,
 $Q_{G}$ is also a Dynkin quiver. Thus,
by \cite [Theorem 7.3]{CH} we have

\medskip

{\bf Corollary 4.6.} Assume that $Q$ is a Dynkin quiver  and
$\lambda=(\lambda_{i})_{i\in I}$ satisfies $\lambda_{i}=\lambda_{j}$
once the vertices $i$ and $j$ are in the same $G$-orbit. 
Then there is a Dynkin quiver $Q'$ such that $\Pi_{Q}^{\lambda}G$
is Morita equivalent to $\Pi_{Q'}$.

\end{document}